\documentclass[12pt]{article}
\usepackage{amssymb,color}
\usepackage{amsmath}
\usepackage{hyperref}
\usepackage{url}

\def\sqr#1#2{{\vcenter{\vbox{\hrule height.#2pt
              \hbox{\vrule width.#2pt height#1pt \kern#1pt \vrule width.#2pt}
              \hrule height.#2pt}}}}
\def\signed #1{{\unskip\nobreak\hfil\penalty50
              \hskip2em\hbox{}\nobreak\hfil#1
              \parfillskip=0pt \finalhyphendemerits=0 \par}}
\def\endpf{\signed {$\sqr69$}}

\def\dbR{{\mathop{\rm l\negthinspace R}}}
\def\dbC{{\mathop{\rm l\negthinspace\negthinspace\negthinspace C}}}

\def \Ker{\mbox{\rm Ker}}
\def\3n{\negthinspace \negthinspace \negthinspace }
\def\2n{\negthinspace \negthinspace }
\def\1n{\negthinspace }

\def\dbC{{\mathbb{C}}}

\def\dbR{{\mathbb{R}}}

\def\dbZ{{\mathbb{Z}}}
\def\={\buildrel \triangle \over =}

\def\a{\alpha}

\def\d{\delta}
\def\e{\varepsilon}

\def\k{\kappa}
\def\l{\lambda}

\def\f{\varphi}

\def\o{\omega}

\def\ns{\noalign{\ss} }

\def\pa{\partial}
\def\R{\mathbb{R}}

\def\cB{{\cal B}}

\def\cE{{\cal E}}
\def\cF{{\cal F}}
\def\cG{{\cal G}}

\def\cJ{{\cal J}}
\def\cK{{\cal K}}
\def\cL{{\cal L}}

\def\cN{{\cal N}}

\def\no{\noindent}

\def\ms{\medskip}
\def\bs{\bigskip}
\def\q{\quad}
\def\qq{\qquad}

\def\max{\mathop{\rm max}}
\def\min{\mathop{\rm min}}

\def\pa{\partial}

\def\wt{\widetilde}
\def\cd{\cdot}
\def\cds{\cdots}

\def\spa{\hbox{\rm span$\,$}}

\def\wh{\widehat}

\def\|{\Big |}
\def\({\Big (}
\def\){\Big )}
\def\[{\Big[}
\def\]{\Big]}
\def\be{\begin{equation}}
\def\bel{\begin{equation}\label}
\def\ee{\end{equation}}
\def\bt{\begin{theorem}}
\def\bcd{\begin{condition}}
\def\ecd{\end{condition}}
\def\et{\end{theorem}}
\def\bc{\begin{corollary}}
\def\ec{\end{corollary}}
\def\bde{\begin{definition}}
\def\ede{\end{definition}}
\def\bl{\begin{lemma}}
\def\el{\end{lemma}}
\def\bp{\begin{proposition}}
\def\ep{\end{proposition}}
\def\br{\begin{remark}}
\def\er{\end{remark}}
\def\ba{\begin{array}}
\def\ea{\end{array}}
\def\ed{\end{document}}
\def\ns{\noalign{\ms}}
\def\ds{\displaystyle}

\def\square#1{\vbox{\hrule\hbox{\vrule height#1%
     \kern#1\vrule}\hrule}}
\def\rectangle#1#2{\vbox{\hrule\hbox{\vrule height#1%
     \kern#2\vrule}\hrule}}

\font\tenbb=msbm10 \font\sevenbb=msbm7
\font\fivebb=msbm5

\newfam\bbfam
\scriptscriptfont\bbfam=\fivebb
\textfont\bbfam=\tenbb
\scriptfont\bbfam=\sevenbb

\newtheorem{lemma}{Lemma}[section]
\newtheorem{remark}{Remark}[section]

\newtheorem{theorem}{Theorem}[section]
\newtheorem{corollary}{Corollary}[section]

\newtheorem{definition}{Definition}[section]
\newtheorem{proposition}{Proposition}[section]
\newtheorem{condition}{Condition}[section]
\makeatletter
   
   \@addtoreset{equation}{section}
\makeatother

\begin{document}

\title{\bf Fredholm Transform and Local Rapid
Stabilization for a Kuramoto-Sivashinsky Equation}

\author{Jean-Michel Coron\thanks{Sorbonne Universit\'{e}s, UPMC Univ Paris 06, UMR 7598, Laboratoire
Jacques-Louis Lions, F-75005, Paris, France
E-mail: \texttt{coron@ann.jussieu.fr}. JMC
was supported by ERC advanced grant 266907
(CPDENL) of the 7th Research Framework
Programme (FP7).} \q and \q Qi
L\"{u}\thanks{School of Mathematics,
Sichuan University, Chengdu 610064, China.
E-mail: \texttt{luqi59@163.com}. QL was
supported by ERC advanced grant 266907
(CPDENL) of the 7th Research Framework
Programme (FP7), the NSF of China under
grant 11101070, the project MTM2011-29306
of the Spanish Science and Innovation
Ministry and the Fundamental Research Funds
for the Central Universities in China under
grants ZYGX2012J115.}}

\date{}

\maketitle

\begin{abstract}\no
This paper is devoted to the study of the local
rapid exponential stabilization problem for a
controlled Kuramoto-Sivashinsky equation on a
bounded interval. We build a feedback control
law to force the solution of the closed-loop
system to decay exponentially to zero with
arbitrarily prescribed decay rates, provided
that the initial datum is small enough. Our
approach uses  a method we introduced for the
rapid stabilization of a Korteweg-de Vries
equation. It relies on the construction of a
suitable integral transform and can be applied
to many other equations.
\end{abstract}

\bs

\no{\bf 2010 Mathematics Subject
Classification}. 93D15, 35Q53. \bs

\no{\bf Key Words}. Kuramoto-Sivashinsky
equation, rapid stabilization,  integral
transform.

 \ms

 \ms

%\newpage

%%%%%%%%%%%%%%%%%%%%%%%%%%%%%%%%%%%%%%%%%%%%%%%%%%%%

\section{Introduction }

%%%%%%%%%%%%%%%%%%%%%%%%%%%%%%%%%%%%%%%%%%%%%%%%%%%%

\q\;Consider the following Kuramoto-Sivashinsky
equation:
\begin{equation}\label{csystem1}
\left\{
\begin{array}{ll}\ds
v_t + v_{xxxx} + \l v_{xx} + vv_x = 0
&\mbox{
in } (0,1)\times (0,+\infty),\\
\ns\ds v(t,0)=v(t,1)=0 &\mbox{ on } (0,+\infty),\\
\ns\ds v_{xx}(t,0)= f(t),\,v_{xx}(t,1)=0
&\mbox{ on }
(0,+\infty),\\
\ns\ds v(0,\cd)=v^0(\cd) &\mbox{ in }(0,1),
\end{array}
\right.
\end{equation}
where $\l>0$ and $v^0(\cd)\in L^2(0,1)$.

For $T>0$, let us define \vspace{-0.2cm}
$$
X_T \= C^0([0,T];L^2(0,1))\cap
L^2(0,T;H^2(0,1)\cap H_0^1(0,1)),
$$
\vspace{-0.2cm}
which is endowed with the
norm
$$
|\cd|_{X_T} =
\big(|\cd|^2_{C^0([0,T];L^2(0,1))} +
|\cd|^2_{L^2(0,T;H^2(0,1)\cap
H_0^1(0,1))}\big)^{\frac{1}{2}}.
$$

We first present the following locally
well-posedness result, which is proved in
Appendix \ref {appendix-wellposed} of this paper.

\begin{theorem}\label{well}
Let $F:L^2(0,1)\to \R$ be a continuous
linear map and let $T_0\in (0,+\infty)$.
Then for given $v^0\in L^2(0,1)$,
there exists at most one solution $v\in
X_{T_0}$ of \eqref{csystem1} with
$f(t)=F(v(t,\cd))$.

Moreover, there exist $r_0>0$ and $C_0>0$ such
that, for every $v^0\in L^2(0,1)$ with
\vspace{-0.1cm}
\begin{equation}\label{v0small}
|v^0|_{L^2(0,1)}\leq r_0,
\end{equation}
there exists one
 solution $v\in X_{T_0}$ of \eqref{csystem1}
with $f(t)=F(v(t,\cd))$ and this solution
satisfies
\begin{equation}\label{estimaelinear}
|v|_{X_T}\leq C_0 |v_0|_{L^2(0,1)}
\end{equation}
\end{theorem}

The Kuramoto-Sivashinsky (K-S for short)
equation was first derived in
\cite{1975-Kuramoto} as a model for
Belouzov-Zabotinskii reaction patterns in
reaction-diffusion systems. It describes  a
lot of physical and chemical systems, such
as unstable flame fronts (see
\cite{1980-Sivashinsky} for example),
falling liquid films (see \cite{1986-Chen}
for example) and interfacial instabilities
between two viscous fluids (see
\cite{1985-Hooper} for example). Since the
pioneer works in
\cite{1988-Foias,1985-Nicolaenko}, the
well-posedness and dynamical properties of
the K-S equation were well studied (see
\cite{2006-Bronski,1993-Collet,2005-Giacomelli,1994-Goodman,2002-Zgliczynski}
and the references therein).

We are interested in the stabilization
problems of the K-S equation. There are
many degrees of freedom to choose the
feedback controls. For instance, one can
choose the internal controls, the Dirichlet
boundary control, the Neumann boundary
control, etc., and the feedback law can be
linear or nonlinear with respect to the
output. Such kind of problems were studied
extensively in the literature. Let us
recall some of them.

In \cite{2000-Armaou,2000-Christofides},
the authors studied the global
stabilization for the K-S equation with
periodic boundary conditions. The feedback
controller acts on the whole domain. They
first considered the ordinary differential
equation approximations of the system,
which accurately describe the dominant
dynamics, and then obtained a local
stabilization result through nonlinear
Galerkin's method.

In the above two works, the control acts in
the whole domain.  In \cite{1998-Hu}, the
authors formulated and solved a robust
boundary control problem for the K-S
equation. In \cite{2001-Liu-Krstic}, with
the assumption that $\l<4\pi^2$, the
authors studied the global stabilization
problems of the K-S equation by a nonlinear
boundary feedback control. The control acts
on any two of the four variables $u$,
$u_x$,$u_{xx}$ and $u_{xxx}$ at the
boundary. It was derived using the
combination of spectral analysis and
Lyapunov techniques which guarantees
$L^2$-global exponential stability. It
seems that their method does not work for
$\l>4\pi^2$.

In \cite{2005-Sakthivel}, with the
assumption of the existence of the global
solution to the K-S equation and  $\l$ is
small, the authors studied the robust
global stabilization of the equation
subject to two nonlinear boundary feedback
controls by Lyapunov techniques.

In \cite{2003-Kobayashi}, under the
assumption of the existence of the global
solution and $\l<1$,  the author studied
the adaptive stabilization of the K-S
equation by four nonlinear boundary
feedback controls. The adaptive stabilizer
was constructed by the concept of high-gain
nonlinear output feedback and the
estimation mechanism of the unknown
parameters.

In this paper, we study the local rapid
stabilization problem of \eqref{csystem1}.
For this purpose, we first put a
restriction on $\l$, i.e., we assume that
\vspace{-0.2cm}
\begin{equation}\label{3.2-eq4}
\l\notin \cN\=\left\{j^2\pi^2 +
k^2\pi^2:\,j,k\in\dbZ^+, j\neq k
\right\}.
\end{equation}
In \eqref{3.2-eq4} and in the following, $\dbZ^+$  denotes the set
of positive integers: $\dbZ^+\=\{1,2,3\ldots\}$.
\begin{remark}
The condition \eqref{3.2-eq4} seems
strange. However, it is natural in the
sense that if $\l\in\cN$, then the
linearized system of \eqref{csystem1} at
zero is not approximately controllable. See
Theorem \ref{th app} in
Appendix~\ref{sec-app-contrllability}.

In \cite{2010-Cerpa}, the same phenomenon of
critical values of $\lambda$  appears in the
study of the boundary null controllability of
the K-S equation with other boundary conditions.
Moreover, a pole shifting method is applied in
order to stabilize the noncritical cases.
Furthermore,   in \cite{2011-Cerpa}, the authors
show that when controlling all the boundary data
at one point ($x=0$ or $x=1$), the equation is
always null controllable.

\end{remark}

Let us introduce an integral transform
$K:L^2(0,1) \to L^2(0,1)$ as \vspace{-0.2cm}
$$
(K v)(x) = \int_0^1 k(x,y)v(y)dy,\q \mbox{
for } v\in L^2(0,1).
$$
Here $k$ is the solution to
\begin{equation}\label{system3}
\left\{
\begin{array}{ll}\ds
k_{xxxx} + \l k_{xx}  -  k_{yyyy} - \l
k_{yy} + a k = a \d(x-y) &\mbox{ in }
(0,1)\times (0,1),\\
\ns\ds k(x,0) = k(x,1) =0 &\mbox{ on } (0,1),\\
\ns\ds k_{yy}(x,0) = k_{yy}(x,1)  =0 &\mbox{ on } (0,1),\\
\ns\ds k_{y}(x,0)=0 &\mbox{ on } (0,1),
\\
\ns\ds k(0,y) = k(1,y)=0 &\mbox{ on }
(0,1),\\
\ns\ds k_{xx}(1,y)=0 &\mbox{ on }
(0,1),\\
\end{array}
\right.
\end{equation}
where $a\in\dbR$ and $\d(x-y)$ denotes the
Dirac measure on the diagonal of the square
$[0,1]\times [0,1]$. The definition of a solution to
\eqref{system3} is given in Section~\ref{secconstructionk}.

We assume that \vspace{-0.2cm}
\begin{equation}\label{3.20-eq1}
a\notin \cN_1\=\Big\{ -k^4\pi^4 + \l
k^2\pi^2 +  j^4\pi^4 - \l j^2\pi^2:\,
j,k\in\dbZ^+ \Big\}.
\end{equation}
Under this assumption, we are going to show the following results.
\begin{enumerate}
\item Equation \eqref{system3} has one and only one solution.  The proof of this result is given in  Section~\ref{secconstructionk}.
\item The operator $I-K$ is invertible. The proof of this result is given in Section~\ref{secinvert}.
\item Assume that
\begin{equation}\label{alargeenough}
a> -j^4\pi^4 + \l j^2\pi^2, \q \forall j\in\dbZ^+.
\end{equation}
Let $\nu>0$ be such that
\begin{equation}\label{nuassezpetit}
0<\nu< a+ j^4\pi^4 - \l j^2\pi^2, \q \forall j\in\dbZ^+.
\end{equation}
Then, for every $v^0\in L^2(0,1)$
satisfying $|v^0|_{L^2(0,1)}\leq r$, if we
define the feedback law $F(\cd)$ by
\vspace{-0.2cm}
\begin{equation}\label{deffeedback}
 F( v)= f(t)\=\int_0^1
k_{xx}(0,y) v(t,y)dy,
\end{equation}
then, for the solution $v$ of \eqref{csystem1}, one has
\begin{equation}\label{decreaseI-K}
|(I-K)v(t)|_{L^2(0,1)} \leq
e^{-\nu t}|(I-K)v^0|_{L^2(0,1)},\q \forall t\in[0,+\infty).
\end{equation}
This result is proved in Section~\ref{secproofth}.
\end{enumerate}
\vspace{0.1cm}

 From the above results, we
get the following theorem.
\begin{theorem}\label{th1}
Let us assume that \eqref{3.2-eq4} hold. Let
$a>0$ be such that \eqref{3.20-eq1} and
\eqref{alargeenough} hold and let $\nu>0$ be
such that \eqref{nuassezpetit} holds. Then,
there exist  $r>0$ and $C>0$ such that, for
every $v^0\in L^2(0,1)$ satisfying
$|v^0|_{L^2(0,1)}\leq r$, the solution $v$ to
\eqref{csystem1} with $f(t)\= F(v(t,\cdot))$ and
$v(0)=v^0$  is defined on $[0,+\infty)$ and
satisfies
\begin{equation}\label{th1-eq1}
|v(t,\cd)|_{L^2(0,1)}\leq Ce^{-\nu
t}|v(0,\cd)|_{L^2(0,1)}, \q\mbox{ for every }\,t
\geq 0.
\end{equation}
\end{theorem}
\begin{remark}
Clearly, for every $\nu>0$, there exists $a>0$
such that \eqref{3.20-eq1} and
\eqref{alargeenough}. Hence we got the rapid
stabilization of our K-S control system
\eqref{csystem1}, i.e., for every $\nu>0$, there
exists a (linear) feedback law such that  $0\in
L^2(0,1)$ is exponential stable for the closed
loop system with an exponential decay rate at
least equal to $\nu$.
\end{remark}
\begin{remark}
There are three differences between  system
\eqref{csystem1} and  systems studied in
\cite{2001-Liu-Krstic,2005-Sakthivel}. The
first one  is that we only employ one
control. The second one is   the feedback
control law is linear. The third one is  we
do not assume that $\l$ is small.

Our method can also be applied to deal with
the rapid stabilization problem of the K-S
equation with other type boundary
conditions and controls, such as
$$
\begin{array}{ll}\ds
y(t,0)=y(t,1)=0,\; y_x(t,0)=f(t),\;
y_x(t,1)=0 \mbox{ for } t\in (0,+\infty),\\
\ns\ds y_{xx}(t,0)=y_{xx}(t,1)=0,\;
y_{xxx}(t,0)=f(t),\; y_{xxx}(t,1)=0 \mbox{
for } t\in (0,+\infty).
\end{array}
$$
The control can also be put on the right
end point of the boundary. In these cases,
one has to modify the set $\cN$ given in
\eqref{3.2-eq4} according to the boundary
conditions.
 The proofs
of the stabilization result for these cases
are quite similar to that given in this
paper and so are omitted.
\end{remark}

The way of constructing the feedback
control in the form of \eqref{deffeedback}
was first introduced in
\cite{2014-Coron-Lu} for obtaining the
local rapid stabilization result for KdV
equations. It was motivated by the
backstepping method (see \cite{KS1} for a
systematic introduction of this method). Some ingredients
of the proofs given here for the K-S equation are also inspired from \cite{2014-Coron-Lu}.

Let us briefly explain the idea for introducing
the transform $K$. For simplicity, let us forget
the nonlinearity in \eqref{csystem1} and
therefore consider the following linear
equation:
\begin{equation}\label{3.2-eq1}
\left\{
\begin{array}{ll}\ds
u_t + u_{xxxx} + \l u_{xx} = 0 &\mbox{
in } (0,1)\times (0,+\infty),\\
\ns\ds u(t,0)=u(t,1)=0 &\mbox{ on } (0,+\infty),\\
\ns\ds u_{xx}(t,0)= g(t),\,u_{xx}(t,1)=0
&\mbox{ on }
(0,+\infty),\\
\ns\ds u(0,\cd)=u^0(\cdot) &\mbox{ in }(0,1).
\end{array}
\right.
\end{equation}
Let \vspace{-0.2cm}
$$
\tilde u(t,x)=u(t,x)-\int_0^1
k(x,y)u(t,y)dy,\q g(t)=\int_0^1
k_{xx}(0,y)u(t,y)dy,
$$
where $k(\cd,\cd)$ is a solution to
\eqref{system3}. If $k(\cd,\cd)$ is smooth
enough, one can easily see that $\tilde
u(\cd,\cd)$ is the solution to
\begin{equation}\label{3.2-eq2}
\left\{
\begin{array}{ll}\ds
\tilde u_t + \tilde u_{xxxx} + \l \tilde
u_{xx} + a\tilde u = 0 &\mbox{
in } (0,1)\times (0,+\infty),\\
\ns\ds \tilde u(t,0)=\tilde u(t,1)=0 &\mbox{ on } (0,+\infty),\\
\ns\ds \tilde u_{xx}(t,0)= 0,\,\tilde
u_{xx}(t,1)=0 &\mbox{ on }
(0,+\infty),\\
\ns\ds \tilde u(0,\cd)=u^0(\cd) - \int_0^1
k(\cd,y)u^0(y)dy &\mbox{ in }(0,1).
\end{array}
\right.
\end{equation}
Let us point out that $|\tilde u|_{L^2(0,1)}$
decays exponentially if we take $a>0$ large
enough and that the exponential decay rate goes
to $+\infty$ as  $a\rightarrow +\infty$.
% ??? give lower bound on the decay rate in a remark later on
To deal with the nonlinear term $vv_x$ requires
further arguments which will be given in
Section~\ref{secproofth}.

It seems that it is very difficult to apply
the backstepping method to solve our
problem. Indeed, let us consider the linear
system \eqref{3.2-eq1} again. If we follow
the backstepping method and choose the
feedback control as
\begin{equation}\label{deffeedback1}
\hat u (t,x)= u(t,x) - \int_x^1 \hat k(x,y)
u(t,y)dy,\ g(t) = \int_0^1 \hat
k_{xx}(0,y)u(t,y)dy,
\end{equation}
then, to guarantee that $\hat u(\cd,\cd)$
is a solution to \eqref{3.2-eq2},   $\hat
k$ should solve
\begin{equation}\label{3.1-eq1}
\left\{
\begin{array}{ll}\ds
\hat k_{xxxx} + \l \hat k_{xx}  -  \hat
k_{yyyy} - \l \hat k_{yy} + a \hat k = 0
&\mbox{ in }
(0,1)\times (0,1),\\
\ns\ds  \hat k(x,1) =  \hat k_{yy}(x,1) =0 &\mbox{ on } (0,1),\\
\ns\ds  \hat k(x,x)  =0 &\mbox{ on } (0,1),\\
\ns\ds \hat k_{xx}(x,x) + \hat
k_{xy}(x,x)=0 &\mbox{ on } (0,1),
\\
\ns\ds 4\hat k_{xxx}(x,x)\! +\! 2\hat
k_{yyy}(x,x)\! +\! 6\hat k_{xxy}(x,x)\! +\!
4\hat k_{xyy}(x,x)\! +\! (\l\!\!-\!\!1)\hat
k_y(x,x)\!=\!0 &\mbox{ on } (0,1).
\end{array}
\right.
\end{equation}
With these five boundary restrictions, the
fourth order equation \eqref{3.1-eq1}
becomes overdetermined. Therefore, it is
not clear whether such a function $\hat
k(\cd,\cd)$ exists.

The rest of this paper is organized as
follows. In Section \ref{secconstructionk},
we establish the well-posedness of equation
\eqref{system3}. In Section
\ref{secinvert}, we show that $I-K$ is an
invertible operator. At last, in Section
\ref{secproofth}, we prove Theorem
\ref{th1}.

\section{Well-posedness of \eqref{system3}}
\label{secconstructionk}

%%%%%%%%%%%%%%%%%%%%%%%%%%%%%%%%%%%%%%%%%%%%%%%%%%%%

\q\;This section is devoted to the study of the
well-posedness of equation \eqref{system3}. We
first introduce the definition of the solution
to \eqref{system3}. Let \vspace{-0.3cm}
\begin{equation}\label{def cE}
\begin{array}{ll}\ds
\cE\=\{\rho \in C^\infty([0,1]\times [0,1]):\, \rho(0,y)=\rho(1,y)=\rho(x,0)=\rho(x,1)=0,\\
\ns\ds \qq \;
\rho_x(0,y)=\rho_{xx}(0,y)=\rho_{xx}(1,y)=\rho_{yy}(x,1)=0\}
\end{array}
\end{equation}
and let $\cG$ be the set of $k\in
H^2((0,1)\times (0,1))\cap
H_0^1((0,1)\times (0,1))$ such that
\begin{gather}
\label{regkx}
  \left( x\in (0,1)\mapsto k_{xx}(x,\cdot)\in L^2(0,1)\right) \in C^0([0,1];L^2(0,1)),
\\
\label{regky}
  \left( y\in (0,1)\mapsto k_{yy}(\cdot,y)\in L^2(0,1)\right) \in C^0([0,1];L^2(0,1)),
\\
\label{ky0}
k_{yy}(\cdot,0)=k_{yy}(\cdot,1)=0 \text{ in
} L^2(0,1).
\end{gather}
We call $k(\cd,\cd)\in \cG$ a solution to
\eqref{system3} if \vspace{-0.2cm}
\begin{equation}
\label{deftransposition}
\begin{array}{ll}\ds
\int_0^1\int_0^1 \big[\rho_{xxxx}(x,y)+\l\rho_{xx}(x,y) - \rho_{yyyy}(x,y) - \l\rho_{yy}(x,y)
+ a \rho(x,y) \big]k(x,y)dxdy\\
\ns\ds - \int_0^1 a \rho(x,x)dx=0, \q\mbox{
for every } \rho\in\cE.
\end{array}
\end{equation}
%

%
%???Should be useless in principle

We have the following well-posedness result
for \eqref{system3}.
\begin{theorem}\label{well-posed lm}
Suppose that \eqref{3.2-eq4} and
\eqref{3.20-eq1} hold. Equation
\eqref{system3} has a unique solution in
$\cG$.
\end{theorem}

Before proceeding our proof, we recall the
 following two results.

\begin{lemma}\cite[Page 45, Theorem 15]{1980-Young-book}\label{lm1}
Let $H$ be a separable Hilbert space and
let $\{e_j\}_{j\in\dbZ^+}$ be an
orthonormal basis for $H$. If
$\{f_j\}_{j\in\dbZ^+}$ is an
$\o$-independent sequence such
that\vspace{-0.1cm}
$$
\sum_{j\in\dbZ^+}|f_j-e_j|_H^2<+\infty,
$$
\vspace{-0.1cm}then $\{f_j\}_{j\in\dbZ^+}$
is a Riesz basis for $H$.
\end{lemma}
\begin{lemma}\cite[Page 40, Theorem 12]{1980-Young-book}\label{lm2}
Let $\{e_j\}_{j\in\dbZ^+}$ be a basis of a
Banach space $X$ and let
$\{f_j\}_{j\in\dbZ^+}$ be the associated
sequence of coefficient functionals. If
$\{b_j\}_{j\in\dbZ^+}$ is complete in $X$
and if
$$
\sum_{j\in\dbZ^+}|e_j-b_j|_{X}|f_j|_{X'}<
+\infty,
$$
then $\{b_j\}_{j\in\dbZ^+}$ is a basis for
$X$ which is equivalent to
$\{e_j\}_{j\in\dbZ^+}$.
\end{lemma}

{\it Proof of Theorem~\ref{well-posed
lm}}\,: We divide the proof  into two
steps:
\begin{itemize}
\item Step 1: proof of the uniqueness of the solution to
\eqref{system3};
\item Step 2: proof of the existence of a solution to \eqref{system3}.
\end{itemize}

\vspace{0.3cm}

\begin{center}
\textbf{Step 1: proof of the uniqueness of
the solution to \eqref{system3}}
\end{center}

Assume that $k_1(\cd,\cd)$ and $k_2(\cd,\cd)$
are two solutions to \eqref{system3}. Let
$k_3(\cd,\cd)=k_1(\cd,\cd)-k_2(\cd,\cd)$. Then
$k_3(\cd,\cd)$ is   such that
\eqref{deftransposition} holds without the last
integral term, i.e., $k_3(\cd,\cd)$ is a
solution to
\begin{equation}\label{system4}
\left\{
\begin{array}{ll}\ds
k_{3,xxxx} + \l k_{3,xx}  -  k_{3,yyyy} -
\l k_{3,yy} + a k_3 = 0 &\mbox{ in }
(0,1)\times (0,1),\\
\ns\ds k_3(x,0) = k_3(x,1) =0 &\mbox{ on } (0,1),\\
\ns\ds k_{3,yy}(x,0) = k_{3,yy}(x,1)  =0 &\mbox{ on } (0,1),\\
\ns\ds k_{3,y}(x,0)=0 &\mbox{ on } (0,1),
\\
\ns\ds k_3(0,y) = k_3(1,y)=0 &\mbox{ on }
(0,1),\\
\ns\ds k_{3,xx}(1,y)=0 &\mbox{ on }
(0,1).\\
\end{array}
\right.
\end{equation}

Let us define an unbounded linear operator
$A: D(A)\subset L^2(0,1)\to L^2(0,1)$ as
follows: \vspace{-0.2cm}
\begin{equation}\label{DA}
\left\{
\begin{array}{ll}\ds
D(A) = \{w\in H^4(0,1):\, w(0)=w(1)=w_{xx}(0)=w_{xx}(1)=0\},\\
\ns\ds A w = -w_{xxxx}-\l
w_{xx},\q\forall\,w\in D(A).
\end{array}
\right.
\end{equation}
$A$ is a self-adjoint operator with compact
resolvent. Simple computations give that the
eigenvalues of $A$ are
\begin{equation}\label{1.22-eq5}
\mu_{j}= -j^4\pi^4 + \l j^2\pi^2,\q
j\in\dbZ^+
\end{equation}
and that the eigenfunction corresponding to the
eigenvalue $\mu_j$ is
\begin{equation}\label{1.17-eq5}
\f_{j}(x)=\sqrt{2}\sin (j\pi x),\q
j\in\dbZ^+.
\end{equation}
Note that, by \eqref{3.2-eq4}, one has
\begin{gather}\label{eigendifferent}
\left(j\in \dbZ^+, \q j\in \dbZ^+ \text{ and } j\not = k\right)\Rightarrow
\left( \mu_j\not =\mu_k\right).
\end{gather}

Let us write
\begin{equation}
\label{k3=} k_3(x,y)=\sum_{j\in\dbZ^+}
\psi_j(x)\f_j(y)
\end{equation}
 for the solution to
\eqref{system4}. Then, $\psi_j$ solves
\begin{equation}\label{system16.1}
\left\{
\begin{array}{ll}\ds
\psi_j'''' + \l\psi_j'' + a \psi_j +  \mu_j
\psi_j = 0
&\mbox{ in }(0,1),\\
\ns\ds \psi_j(0) =
\psi_j(1)=\psi_{j}''(1)=0.
\end{array}
\right.
\end{equation}
Let $c_j \=  \psi_{j}''(0)$ ($j\in
\mathbb{Z}^+$). We consider the following
equation:
\begin{equation}\label{system16}
\left\{
\begin{array}{ll}\ds
\check\psi_j'''' + \lambda \check\psi_j'' + a
\check\psi_j + \mu_j \check\psi_j = 0
&\mbox{ in }(0,1),\\
\ns\ds \check\psi_j(0) = \check\psi_j(1)=\check\psi_j''(1)=0,\\
\ns\ds  \check\psi_j''(0)=1.
\end{array}
\right.
\end{equation}
Since, by \eqref{3.20-eq1} and \eqref{1.22-eq5}, $a+\mu_j$ is not an eigenvalue of $A$,
\eqref{system16} has a unique solution.
Moreover, $\psi_j=c_j\check\psi_j$  for
every $j\in\dbZ^+$. The four roots of
\vspace{-0.2cm}
$$
r^4 + \l r^2 + a + \mu_{j}=0
$$
are \vspace{-0.4cm}
\begin{equation}\label{1.22-eq13}
\begin{array}{ll}\ds
r^{(1)}_{j}=\(\frac{-\l+\sqrt{\l^2-4(a+\mu_{j})}}{2}\)^{\frac{1}{2}},
\;r^{(2)}_{j}=-\(\frac{-\l+\sqrt{\l^2-4(a+\mu_{j})}}{2}\)^{\frac{1}{2}}=-r^{(1)}_j,\\
\ns\ds r^{(3)}_{j} =
\(\frac{-\l-\sqrt{\l^2-4(a+\mu_{k,j})}}{2}\)^{\frac{1}{2}},\;
r^{(4)}_{j}=-\(\frac{-\l-\sqrt{\l^2-4(a+\mu_{j})}}{2}\)^{\frac{1}{2}}=-r^{(3)}_j.
\end{array}
\end{equation}
Easy computations show that there exists $C>0$
such that, with $i\=\sqrt{-1}$,
\begin{equation}\label{6.17-eq2}
\left|r^{(1)}_{j} - j\pi +
\frac{\l}{2j\pi}\right|\leq \frac{C}{j^3},\q
\left|ir^{(4)}_{j} -
\left(j\pi-\frac{a}{4j^3\pi^3}\right)\right|\leq
\frac{C}{j^5},  \q \forall j\in j\in\dbZ^+
\end{equation}
Let us assume that
\begin{equation}\label{qssumptiopourdiff}
a+\mu_j\not =0 \text{ and } \lambda^2\not = 4(a+\mu_j).
\end{equation}
The cases where \eqref{qssumptiopourdiff} does
not hold require slight modifications in the
arguments given below. We omit them. Note that
\eqref{qssumptiopourdiff} holds for $j$ large
enough. From \eqref{1.22-eq13} and
\eqref{qssumptiopourdiff}, one gets that the
four roots $r^{(1)}_{j}$, $r^{(2)}_{j}$,
$r^{(3)}_{j}$ and $r^{(4)}_{j}$ are distinct.
Hence there exists four complex numbers
$\a_j^{(1)}$, $\a_j^{(2)}$, $\a_j^{(3)}$ and
$\a_j^{(4)}$ such that  \vspace{-0.2cm}
\begin{equation}\label{eqcheckpsi}
\check\psi_j = \a_j^{(1)}\cosh(r_j^{(1)}x) +
\a_j^{(2)}\sinh(r_j^{(1)}x) +
\a_j^{(3)}\cos(ir_j^{(3)}x) +
\a_j^{(4)}\sin(ir_j^{(4)}x).
\end{equation}
Let us emphasize that throughout this section,
the functions as well as the sequences of
numbers are complex valued. However, at the end
of this section we will check that $k$ is real
valued. From the boundary conditions of
\eqref{system16} and \eqref{eqcheckpsi}, we get
that
\begin{equation}\label{1.22-eq22}
\left\{
\begin{array}{ll}\ds
\a_j^{(1)} + \a_j^{(3)}=0,\\
\ns\ds \a_j^{(1)}\cosh(r_j^{(1)}) +
\a_j^{(2)}\sinh(r_j^{(1)}) +
\a_j^{(3)}\cos(r_j^{(3)}) +
\a_j^{(4)}\sin(r_j^{(4)})=0,\\
\ns\ds   \a_j^{(1)}(r_j^{(1)})^2
+ \a_j^{(3)}(r_j^{(4)})^2 =1,\\
\ns\ds \a_j^{(1)}(r_j^{(1)})^2\cosh(r_j^{(1)}) +
\a_j^{(2)}(r_j^{(1)})^2 \sinh(r_j^{(1)}) -
\a_j^{(3)}(r_j^{(4)})^2\cos(ir_j^{(4)})\\
\ns\ds -
\a_j^{(4)}(r_j^{(4)})^2\sin(ir_j^{(4)})=0.
\end{array}
\right.
\end{equation}
By means of the first equation of
\eqref{1.22-eq22}, we find that $\a_j^{(1)} =-
\a_j^{(3)}$. This, together with
\eqref{1.22-eq13} and the third equation of
\eqref{1.22-eq22}, implies that \vspace{-0.2cm}
\begin{equation}\label{2.10-eq1}
\a_j^{(1)}=-
\a_j^{(3)}=\(\sqrt{\l^2-4(a+\mu_j)}\)^{-1}.
\end{equation}
(Note that by \eqref{qssumptiopourdiff}, the
right hand side of \eqref{2.10-eq1} is well
defined.) From \eqref{qssumptiopourdiff} and
\eqref{qssumptiopourdiff}, one gets
\begin{equation}\label{detnonzeroalpha}
(r_j^{(1)})^2+(r_j^{(4)})^2\not =0
\end{equation}
According to  \eqref{1.22-eq13},
\eqref{2.10-eq1}, the second and fourth
equations in \eqref{1.22-eq22}  and
\eqref{detnonzeroalpha}, we find that
\begin{equation}\label{2.10-eq2}
\a_j^{(2)}=-\a_j^{(1)}\coth(r_j^{(1)}),\q
\a_j^{(4)}=-\a_j^{(3)}\cot(ir_j^{(4)}).
\end{equation}
Thus, we obtain that \vspace{-0.3cm}
\begin{equation}\label{2.10-eq3}
\begin{array}{ll}\ds
\check\psi_j \3n&\ds =
\a_j^{(1)}\cosh(r_j^{(1)}x) - \a_j^{(1)}
\coth(r_j^{(1)})\sinh(r_j^{(1)}x)
\\
\ns&\ds \q - \a_j^{(1)}\cos(ir_j^{(4)}x) +
\a_j^{(1)}\cot(ir_j^{(4)})\sin(ir_j^{(4)}x).
\end{array}
\end{equation}
Note that, by \eqref{3.20-eq1},
\begin{equation}\label{anot0}
  a\not =0.
\end{equation}
Due to  \eqref{6.17-eq2} and \eqref{anot0}, we
have \vspace{-0.2cm}
\begin{equation}\label{2.10-eq4}
\cot(ir_j^{(4)})=O(j^3) \mbox{ as }
j\to\infty.
\end{equation}

Let \vspace{-0.2cm}
\begin{equation}\label{2.10-eq5}
\hat\psi_j =
\frac{\sqrt{2}}{\a_j^{(1)}\cot(ir_j^{(4)})}\check\psi_j
\mbox{ for } j\in\dbZ^+.
\end{equation}
Let us now prove the following lemma.

\begin{lemma}
\label{lemRiesz}
The family $\{\hat\psi_j\}_{j\in\dbZ^+ }$ is a Riesz
basis of $L^2(0,1)$.
\end{lemma}
{\it Proof of Lemma~\ref{lemRiesz}.} First, we claim that
\begin{equation}\label{4.18-eq4}
\sum_{j\in\dbZ^+} \big|\hat\psi_{j} - \f_j
\big|_{L^2(0,1)}^2 <\infty.
\end{equation}
Indeed, from  \eqref{6.17-eq2} and
\eqref{2.10-eq3} to \eqref{2.10-eq5}, we see
that there exists $C >0$ such that, for every $j
\in\dbZ^+$,
\begin{equation}\label{5.15-eq5}
\begin{array}{ll}\ds
\q \int_{0}^1 |\hat \psi_{j} -  \f_j|^2 dx \\
\ns\ds \leq 2\( \int_{0}^1 \big|
\sqrt{2}\sin(ir_j^{(4)}x) -
\sqrt{2}\sin(j\pi x)\big|^2 dx \\
\ns\ds \q + \int_{0}^1
\frac{2}{|\cot(ir_j^{(4)})|^2}
\big|\cosh(r_j^{(1)}x) -
\coth (r_j^{(1)})\sinh(r_j^{(1)}x) -
\cos(ir_j^{(4)}x)\big|^2 dx\)
\\
\ns \ds \leq  \frac{C}{j^6}.
\end{array}
\end{equation}
Since $ \ds\sum_{j=1}^{+\infty} j^{-6} <+\infty$, we have
\eqref{4.18-eq4}.

Let $T>0$. Consider the following  control
system:
\begin{equation}\label{10.4-eq1}
\left\{
\begin{array}{ll}\ds
\vartheta_t + \vartheta_{xxxx} + \l\vartheta_{xx} = 0 &\mbox{ in } (0,T)\times(0,1),\\
\ns\ds \vartheta(t,0)=\vartheta(t,1)=0 &\mbox{ in } (0,T),\\
\ns\ds \vartheta_{xx}(t,0)=\eta(t),\;
\vartheta_{xx}(t,1)=0 &\mbox{ in } (0,T),
\end{array}
\right.
\end{equation}
where $\eta(\cd)\in L^2(0,T)$ is the
control. Let $\tilde
\vartheta=A^{-1}\vartheta$ (recall
\eqref{DA} for the definition of $A$).
Then, $\tilde \vartheta$ solves
\vspace{-0.2cm}
\begin{equation}\label{10.4-eq2}
\pa_t \tilde \vartheta = A\tilde \vartheta -
\eta(t)b,
\end{equation}
Here $b(\cd)$ is the solution to
\begin{equation}\label{10.31-eq1}
\left\{
\begin{array}{ll}\ds
b_{xxxx} + \l b_{xx} =0 &\mbox{ in }(0,1),\\
\ns\ds
b(0)=b(1)=0,\;b_{xx}(0)=1,\;b_{xx}(1)=0.
\end{array}
\right.
\end{equation}
Clearly, $b\in L^2(0,1)$. Let
$b=\sum_{j\in\dbZ^+}b_j\f_j$. Since  system
\eqref{10.4-eq1} is approximately
controllable in $L^2(0,1)$ (see Theorem
\ref{th app}), we get that \eqref{10.4-eq1}
is also approximately controllable in
$A^{-1}\big(L^2(0,1)\big)$. In particular,
\begin{equation}\label{bjnot0}
b_j\neq 0, \q j\in\dbZ^+.
\end{equation}

For $j\in\dbZ^+$, since, as already used above,
$a + \mu_j$ is not an eigenvalue of $A$ (recall
once more \eqref{3.20-eq1} and
\eqref{1.22-eq5}), there exists a
$\tilde\psi_j\in D(A)$  such that
\begin{equation}\label{10.4-eq3}
-A\tilde\psi_j + (a + \mu_j)\tilde\psi_j =
-\frac{\sqrt{2}}{\a_j^{(1)}\cot(ir_j^{(4)})}(a+\mu_j)
b.
\end{equation}
Thanks to  \eqref{system16}, \eqref{2.10-eq5}
and \eqref{10.4-eq3}, for every $j\in\dbZ^+$, we
have \vspace{-0.2cm}
\begin{gather}\label{10.4-eq3.1}
A^{-1}\tilde\psi_j =  (a+\mu_j)^{-1}\tilde\psi_j
-
\frac{\sqrt{2}}{\a_j^{(1)}\cot(ir_j^{(4)})}A^{-1}b,
\\
\label{10.4-eq10}
\hat\psi_j = \tilde \psi_j +
\frac{\sqrt{2}}{\a_j^{(1)}\cot(ir_j^{(4)})}b.
\end{gather}
Let us assume that there exists
$\{a_j\}_{j\in\dbZ^+}\in \ell^2(\dbZ^+)$
such that $
\sum_{j\in\dbZ^+}a_j\hat\psi_j=0. $ Then, from
\eqref{10.4-eq10}, we obtain that
\begin{equation}\label{10.4-eq5.1}
\sum_{j\in\dbZ^+}a_j\(\tilde \psi_j +
\frac{\sqrt{2}}{\a_j^{(1)}\cot(ir_j^{(4)})}
b\)=0.
\end{equation}
Applying $A^{-1}$ to \eqref{10.4-eq5.1},
and using \eqref{10.4-eq3.1}, one gets that
\begin{equation}\label{10.4-eq6}
\sum_{j\in\dbZ^+}
a_j(a+\mu_j)^{-1}\tilde\psi_j=0,
\end{equation}
which, together with  \eqref{10.4-eq10}, deduces
that
\begin{equation}\label{10.4-eq6.1}
\[\sum_{j\in\dbZ^+}
a_j\frac{\sqrt{2}}{\a_j^{(1)}\cot(ir_j^{(4)})}(a+\mu_j)^{-1}\]b-\sum_{j\in\dbZ^+}
a_j(a+\mu_j)^{-1}\hat\psi_j=0.
\end{equation}
Applying $A^{-1}$ to \eqref{10.4-eq6} and
using \eqref{10.4-eq3.1} again, we find
that
\begin{equation}\label{10.4-eq7}
\[\sum_{j\in\dbZ^+}
a_j\frac{\sqrt{2}}{\a_j^{(1)}\cot(ir_j^{(4)})}(a+\mu_j)^{-1}\]A^{-1}b
-\sum_{j\in\dbZ^+}
a_j(a+\mu_j)^{-2}\tilde\psi_j=0,
\end{equation}
which, together with \eqref{10.4-eq10}, gives
\begin{equation}\label{10.4-eq7.1}
\begin{array}{ll}\ds
\[\sum_{j\in\dbZ^+}
a_j\frac{\sqrt{2}}{\a_j^{(1)}\cot(ir_j^{(4)})}(a+\mu_j)^{-1}\]A^{-1}b
+  \[\sum_{j\in\dbZ^+}
a_j\frac{\sqrt{2}}{\a_j^{(1)}\cot(ir_j^{(4)})}(a+\mu_j)^{-2}\]
b \\
\ns\ds - \sum_{j\in\dbZ^+}
a_j(a+\mu_j)^{-2}\hat\psi_j=0.
\end{array}
\end{equation}
By induction, one gets that, for every
positive integer $p$,
\begin{equation}\label{eqp}
\begin{array}{ll}\ds
\(\sum_{j\in\dbZ^+}
a_j\frac{\sqrt{2}}{\a_j^{(1)}\cot(ir_j^{(4)})}(a+\mu_j)^{-1}
\)A^{-p}b \\
\ns\ds + \sum_{k=2}^{p}\(\sum_{j\in\dbZ^+}
a_j\frac{\sqrt{2}}{\a_j^{(1)}\cot(ir_j^{(4)})}(a+\mu_j)^{-k}
\)A^{k-p-1}b \\
\ns\ds +   \sum_{j\in\dbZ^+}
a_j\frac{\sqrt{2}}{\a_j^{(1)}\cot(ir_j^{(4)})}(a+\mu_j)^{-p-1}b
-
 \sum_{j\in\dbZ^+}
a_j(a+\mu_j)^{-p-1} \hat\psi_j=0.
\end{array}
\end{equation}
If \vspace{-0.3cm}
\begin{equation}
\label{premieresomme} \sum_{j\in\dbZ^+}
a_j\frac{\sqrt{2}}{\a_j^{(1)}\cot(ir_j^{(4)})}(a+\mu_j)^{-1}\not
=0,
\end{equation}
we get from \eqref{eqp} that
\vspace{-0.4cm}
\begin{equation}
\label{inclusionpsi}
\{A^{-p}b\}_{p\in(\{0\}\cup\dbZ^+)}\subset\spa\{\hat\psi_j\}_{j\in\dbZ^+}.
\end{equation}
If
$\overline{\spa\{\hat\psi_j\}_{j\in\dbZ^+}}\neq
L^2(0,1)$, then we can find a nonzero
function $ d=\sum_{j\in\dbZ^+}d_j\f_j\in
L^2(0,1)$ such that \vspace{-0.2cm}
\begin{equation}
\label{scalarproduct0}
(h,d)_{L^2(0,1)}=0,\;\mbox{ for every }h\in
\spa\{\hat\psi_j\}_{j\in\dbZ^+}.
\end{equation}
 From \eqref{inclusionpsi} and \eqref{scalarproduct0}, we obtain that $(A^{-p}b,d)_{L^2(0,1)}=0$
for every $p\in(\{0\}\cup\dbZ^+)$.
Therefore, we get \vspace{-0.2cm}
\begin{equation}\label{10.4-eq8}
\sum_{j\in\dbZ^+}b_j\mu_j^{-p}d_j=0 \mbox{
for all } p\in(\{0\}\cup\dbZ^+).
\end{equation}
Let us define a complex variable function
$G(\cd):\mathbb{C}\to \mathbb{C}$ as
\vspace{-0.2cm}
$$
G(z)=\sum_{j\in\dbZ^+}d_jb_je^{\mu_j^{-1}z},\q
z\in\dbC.
$$
Then, it is clear that $G(\cd)$ is a
holomorphic function. From
\eqref{10.4-eq8}, we see that $
G^{(p)}(0)=0 \;\mbox{ for every
}p\in(\{0\}\cup\dbZ^+) $. Thus, we obtain
that
\begin{equation}\label{Gnull}
G(\cd)= 0.
\end{equation}
Using \eqref{Gnull} and \eqref{eigendifferent},
and looking at the asymptotic behavior of $G(z)$
as $z\in \R$ tends to $+\infty$ -if some $
\mu_j$ are positive- and to $-\infty$ one gets
that $d_jb_j=0$ for all $j\in\dbZ^+$. Since
$b_j\neq 0$, we  get that $d_j=0$ for all
$j\in\dbZ^+$. Therefore, we get $d=0$, which
leads to a contradiction. Hence,
\eqref{premieresomme} implies that
\begin{equation}
\label{density}
\overline{\spa\{\hat\psi_j\}_{j\in\dbZ^+}}=L^2(0,1).
\end{equation}
If \eqref{premieresomme} does not hold but
\vspace{-0.4cm}
$$
\sum_{j\in\dbZ^+}
a_j\frac{\sqrt{2}}{\a_j^{(1)}\cot(ir_j^{(4)})}(a+\mu_j)^{-2}\neq
0,
$$
then by using \eqref{eqp} again, we obtain
that \vspace{-0.2cm}
$$
\{A^{-p}b\}_{p\in
(\{0\}\cup\dbZ^+)}\subset\spa\{\hat\psi_j\}_{j\in\dbZ^+}.
$$
By a similar argument, we find that
\eqref{density} again holds. Similarly, we
can get that, if there is a $p\in\dbZ^+$
such that \vspace{-0.2cm}
\begin{equation}\sum_{j\in\dbZ^+} a_j
\frac{\sqrt{2}}{\a_j^{(1)}\cot(ir_j^{(4)})}(a+\mu_j)^{-p}
\neq 0,
\end{equation}
then \eqref{density} holds.

On the other hand, if \vspace{-0.2cm}
$$
\sum_{j\in\dbZ^+} a_j
\frac{\sqrt{2}}{\a_j^{(1)}\cot(ir_j^{(4)})}(a+\mu_j)^{-p}
= 0 \mbox{ for every
}p\in(\{0\}\cup\dbZ^+),
$$
we  define a function \vspace{-0.2cm}
$$
\wt
G(z)\=\sum_{j\in\dbZ^+}a_j\frac{\sqrt{2}}{\a_j^{(1)}\cot(ir_j^{(4)})}(a+\mu_j)^{-1}
e^{(a+\mu_j)^{-1}z} \mbox{ for every
}z\in\dbC,
$$
and it is clear that $\wt G(\cd)$ is a
holomorphic function and $ \wt
G^{(p)}(0)=0$ for every
$p\in\{0\}\cup\dbZ^+$,
which implies that $\wt G(\cd)= 0$. Therefore,
as above, we conclude that $a_j=0$ for every
$j\in\dbZ^+$.

By the above argument, we know that either
$\{\hat\psi_j\}_{j\in\dbZ^+}$ is
$\o$-independent or it is complete in
$L^2(0,1)$.

We first deal with the case that
$\{\hat\psi_j\}_{j\in\dbZ^+}$ is
$\o$-independent. Let us take $H=L^2(0,1)$
and put $e_{j} = \f_j$, $f_{j} =
\hat\psi_{j}$ for $j\in\dbZ^+$ in Lemma
\ref{lm1}. Then, by \eqref{5.15-eq5}, the
conditions of Lemma \ref{lm1} are
fulfilled. Thus,
$\{\hat\psi_j\}_{j\in\dbZ^+}$ is a Riesz
basis of $L^2(0,1)$. This concludes the proof of Lemma~\ref{lemRiesz}.

Next, we consider the case that
$\{\hat\psi_j\}_{j\in\dbZ^+}$ is complete
in $L^2(0,1)$. In Lemma
\ref{lm2}, let us set $X\=L^2(0,1)$, $
e_{j} \= \f_j$, $f_{j} \= \f_j$, $b_{j} \=
\hat\psi_{j}$ for $j\in\dbZ^+$.  Then, by \eqref{5.15-eq5}, it is
easy to see that the conditions of Lemma
\ref{lm2} are fulfilled. Therefore,
 $\{\hat\psi_j\}_{j\in\dbZ^+}$ is a Riesz
basis of $L^2(0,1)$.

Now we  estimate $\{c_j\}_{j\in\dbZ^+}$.
From \eqref{1.17-eq5}, we get that $
\int_0^1 |\f_{j,yy}(y)|^2 dy =j^4\pi^4$.
This, together with the fact that
$\{\hat\psi_j\}_{j\in\dbZ^+}$ is a Riesz
basis of $L^2(0,1)$ and $k_3(\cd,\cd)\in
H^2((0,1)\times (0,1))$, implies that
$$
\begin{array}{ll}\ds
+\infty\3n&\ds> \int_0^1\int_0^1
|k_{3,yy}(x,y)|^2dxdy \geq \int_0^1\int_0^1
\|\sum_{j\in\dbZ^+}
c_j\frac{\sqrt{2}}{\a_j^{(1)}\cot(ir_j^{(4)})}\hat\psi_j(x)\f_{j}''(y)\|^2dxdy\\
\ns&\ds \geq C\int_0^1 \sum_{j\in\dbZ^+}\|
c_j\frac{\sqrt{2}}{\a_j^{(1)}\cot(ir_j^{(4)})}
\f_{j}''(y)\|^2 dy\geq C\sum_{j\in\dbZ^+}
\Big|c_j\frac{\sqrt{2}}{\a_j^{(1)}\cot(ir_j^{(4)})}j^2\Big|^2.
\end{array}
$$
Hence, we find that \vspace{-0.2cm}
\begin{equation}\label{10.5-eq7}
\bigg\{c_j\frac{\sqrt{2}}{\a_j^{(1)}\cot(ir_j^{(4)})}j^2\bigg\}_{j\in\dbZ^+}\in
\ell^2(\dbZ^+).
\end{equation}
From \eqref{1.17-eq5},  we get  that
\begin{equation}
\label{asympfjy-better} \f_{j}'(0)=\sqrt{2}j\pi,
\q \forall j\in \mathbb{Z}^+.
\end{equation}
From \eqref{10.5-eq7} and
\eqref{asympfjy-better}, we obtain that
\begin{equation}\label{10.5-eq8}
\bigg\{c_j\frac{\sqrt{2}}{\a_j^{(1)}\cot(ir_j^{(4)})}\f_{j}'(0)\bigg\}_{j\in\dbZ^+}\in
\ell^2(\dbZ^+).
\end{equation}
Using \eqref{k3=}, $k_{3,y}(x,0)=0$ and the
second inequality of
\eqref{asympfjy-better}, we find that
\vspace{-0.2cm}
$$
\sum_{j\in\dbZ^+}
\psi_j(\cd)\f_{j}'(0)=\sum_{j\in\dbZ^+}
c_j\frac{\sqrt{2}}{\a_j^{(1)}\cot(ir_j^{(4)})}\f_{j}'(0)\hat\psi_j(\cd)=0
\q\mbox{ in } L^2(0,1).
$$
Since $\{\hat\psi_j\}_{j\in \dbZ }$ is a
Riesz basis of $L^2(0,1)$, we get that
\vspace{-0.2cm}
\begin{equation}\label{cjfj}
c_j\frac{\sqrt{2}}{\a_j^{(1)}\cot(ir_j^{(4)})}\f_{j}'(0)=0
\mbox{ for every } j\in\dbZ^+.
\end{equation}
From \eqref{asympfjy-better}, we have
$\f_{j}'(0)\not=0$, which together with
\eqref{cjfj},
 gives  that $c_j=0$ for every $j\in\dbZ^+$. This
implies that equation \eqref{system4}
admits a unique solution $k_3=0$.
Therefore, we obtain that equation
\eqref{system3} admits at most one
solution. This concludes Step 1.

\vspace{0.3cm}

\begin{center}
\textbf{Step 2: proof of the existence of a
solution to \eqref{system3}}
\end{center}

Denote by $D(A)'$ the dual space of $D(A)$
with respect to the pivot space $L^2(0,1)$.
Let
$h(\cd)=\sum_{j\in\dbZ^+}h_j\f_j(\cd)\in
D(A)$, i.e., $|h(\cd)|_{D(A)}^2 =
\sum_{j\in\dbZ^+}|h_j\mu_j|^2<+\infty$. We
have that \vspace{-0.25cm}
$$
\sum_{j\in\dbZ^+}|h_j\f_{j}'(0)| \leq
\(\sum_{j\in\dbZ^+}|h_j\mu_j|^2\)^\frac{1}{2}\(\sum_{j\in\dbZ^+}|\f_{j}'(0)\mu_j^{-1}|^2\)^\frac{1}{2}<+\infty.
$$
Hence\vspace{-0.35cm}
\begin{equation}
\sum_{j=1}^{n}\f_{j}'(0) \f_j(\cd) \text{ is
convergent in $D(A)'$ as  $n$ tends to
$+\infty$},
\end{equation}
which  allows us to define
\begin{equation}
\label{defhatb} \hat b(\cd)
\=\sum_{j\in\dbZ^+}\f_{j}'(0) \f_j(\cd) \in
D(A)'.
\end{equation}
Furthermore, it is clear that \vspace{-0.35cm}
$$
(h,\hat b)_{D(A),D(A)'} =
\sum_{j\in\dbZ^+}h_j\f_{j}'(0)=h'(0)= -\d'_0(h),
$$
\vspace{-0.25cm}which implies that
\begin{equation}\label{11.8-eq5}
\hat b = -\d'_0 \q\mbox{ in } D(A)'.
\end{equation}
Let
\begin{equation}\label{11.7-eq1}
a_{j} \= \frac{a}{\f'_{j}(0)},\;\;\phi_j \=
\f_{j} + (-A+\mu_{j}+a)^{-1}(a_{j}\hat b - a
 \f_{j}) \in L^2(0,1) \;\mbox{ for } j\in\dbZ^+.
\end{equation}
From  \eqref{1.17-eq5} and
\eqref{11.7-eq1}, we have that
\vspace{-0.2cm}
\begin{equation}\label{10.30-eq1}
\begin{array}{ll}\ds
\q\int_{0}^1 \Big|\f_{j}(x) -
\phi_{j}(x)\Big|^2 dx\\
\ns\ds = \int_{0}^1
\big|(-A+\mu_{j}+a)^{-1}(a_{j}\hat b -
a \f_{j})\big|^2 dx\\
\ns\ds =  \int_{0}^1
\Big|\sum_{k\in\dbZ^+\setminus\{j\}}
(-\mu_k+\mu_{j}+a)^{-1}\frac{a\f_k'(0)}{\f_{j}'(0)}
\f_k(x)\Big|^2dx
\\
\ns\ds = \sum_{k\in\dbZ^+\setminus\{j\}}
\Big|(-\mu_k+\mu_{j}+a)^{-1}\frac{a\f_k'(0)}{\f_{j}'(0)} \Big|^2\\
\ns\ds\leq  C \sum_{k\in\dbZ^+\setminus\{j\}}
\left|k^4-j^4 \right|^{-2}j^{-2}k^2 \q\mbox{ for
all }j\in\dbZ^+.
\end{array}
\end{equation}
From \eqref{10.30-eq1}, for
$j>0$,\vspace{-0.1cm}
\begin{equation}\label{11.7-eq2}
\begin{array}{ll}\ds
\q\int_{0}^1 \Big|\f_{j}(x) -
\phi_{j}(x)\Big|^2 dx
\\
\ns\ds \leq
C\bigg(\sum_{k>2j}\frac{k^2}{\big|k^4-j^4\big|^{2}j^{2}}
+ \sum_{j< k\leq
2j}\frac{k^2}{\big|k^4-j^4\big|^{2}j^{2}} +
\sum_{0< k
<j}\frac{k^2}{\big|k^4-j^4\big|^{2}j^{2}}\bigg).
\end{array}
\end{equation}
Now we estimate the terms in the right hand
side of \eqref{11.7-eq2}. First, we have
that
\begin{equation}\label{11.7-eq3}
\sum_{k>2j}\frac{k^2}{\big|k^4-j^4\big|^{2}j^{2}}\leq
\sum_{k>2j}\frac{16}{15k^{6}j^{2}} \leq
\frac{16}{15j^6}\sum_{k>2j}\frac{1}{k^{2}}\leq
\frac{1}{j^7}.
\end{equation}
Next,
\begin{equation}\label{11.7-eq4}
\begin{array}{ll}\ds
\q\sum_{j< k\leq
2j}\frac{k^2}{\big|k^4-j^4\big|^{2}j^{2}}\leq
4\sum_{j< k\leq
2j}\frac{1}{\big|k^4-j^4\big|^{2}} =
4\sum_{1\leq l\leq
j}\frac{1}{\big|(j+l)^4-j^4\big|^{2}}\\
\ns\ds \leq 4\sum_{1\leq l\leq
j}\frac{1}{\big|4j^3l+6j^2l^2\big|^{2}}\leq
\frac{1}{4}\sum_{-j\leq l\leq
-1}\frac{1}{j^4}\frac{1}{(j+l)^2l^2}=
\frac{1}{4}\sum_{1\leq l\leq
j}\frac{1}{j^4}\[\frac{1}{j}\(\frac{1}{l}-\frac{1}{j+l}\)\]^2\\
\ns\ds \leq \frac{1}{2}\sum_{1\leq l\leq
j}\frac{1}{j^6}\(\frac{1}{l^2}+\frac{1}{(j-l)^2}\)\leq
\frac{\pi^2}{6}\frac{1}{j^6}.
\end{array}
\end{equation}
Similarly, we can obtain that
\vspace{-0.2cm}
\begin{equation}\label{11.7-eq5}
\sum_{0<k<j}\frac{k^2}{\big|k^4-j^4\big|^{2}j^{2}}
\leq \frac{\pi^2}{6}\frac{1}{j^6}.
\end{equation}
From \eqref{11.7-eq2} to \eqref{11.7-eq5},
we know that there is a constant $C>0$ such
that, for all positive integer $j$,
\begin{equation}\label{10.30-eq7}
\int_{0}^1 \Big|\f_{j}(x) -
\phi_{j}(x)\Big|^2 dx \leq \frac{C}{j^6}.
\end{equation}

Furthermore, by similar arguments, we can obtain
that $\phi_j\in H^2(0,1)$ and that
\begin{equation}\label{10.30-eq2}
\int_{0}^1 \Big|\f_{j}'(x) - \phi_{j}'(x)\Big|^2
dx \leq \frac{C}{j^4}\q\mbox{ for all
}j\in\dbZ^+
\end{equation}
and
\begin{equation}\label{10.30-eq2-double}
\int_{0}^1 \Big|\f_{j}''(x) -
\phi_{j}''(x)\Big|^2 dx \leq
\frac{C}{j^2}\q\mbox{ for all }j\in\dbZ^+.
\end{equation}
Let us check that
\begin{equation}\label{phijH6}
  \phi_j\in H^6(0,L) \text { and } \phi_j''(1)=0, \q \forall j\in \dbZ^+.
\end{equation}
Simple computations show that
\begin{equation}\label{1surksinkx}
\sum_{k=1}^{+\infty}\frac{1}{k}\sin(k\pi x)=\frac{\pi}{4} (1-x) \text{ in } L^2(0,1).
\end{equation}
Let $j\in \dbZ^+$. From \eqref{1.22-eq5}, there
exists $C>0$ (depending on $j$) such that
\begin{equation}\label{1surksinkx1}
\left|\frac{1}{-\mu_k+\mu_j+a}-\frac{1}{\pi^4k^4}\right|\leq \frac{C} {k^8},\q \forall k\in \dbZ^+.
\end{equation}
The two statements of \eqref{phijH6} follow from
\eqref{1.17-eq5}, \eqref{defhatb},
\eqref{11.7-eq1}, \eqref{1surksinkx} and
\eqref{1surksinkx1}.

With the same strategy to prove that
$\{\hat\psi_j\}_{j\in\dbZ^+}$ is a Riesz basis
of $L^2(0,1)$, we  can also show that
\begin{equation}\label{10.31-eq3}
\Big\{-(a+\mu_{j})A^{-1}\phi_{j}\Big\}_{j\in\dbZ^+}
\mbox{ is a Riesz basis   of }L^2(0,1).
\end{equation}
From \eqref{1.22-eq5} and
\eqref{asympfjy-better}, we  get that $
\sum_{j\in\dbZ^+}\f_{j}'(0)\mu_j^{-1}\f_j$ is
convergent in   $L^2(0,1)$. From
\eqref{10.31-eq3}, there is $\{\hat
c_j\}_{j\in\dbZ^+}\in l^2(\dbZ^+)$ such that
\vspace{-0.2cm}
$$
-\sum_{j\in\dbZ^+} \hat c_j
(a+\mu_{j})A^{-1}\phi_{j}=\sum_{j\in\dbZ^+}\f_{j}'(0)\mu_j^{-1}\f_j
\q\mbox{ in }\; L^2(0,1),
$$
which, together with \eqref{defhatb} and
\eqref{11.8-eq5}, implies that
\begin{equation}\label{11.2-eq2}
-\sum_{j\in\dbZ^+} \hat c_j (a+\mu_{j})
\phi_{j}=\sum_{j\in\dbZ^+}\f_{j}'(0)\f_j=-\d'_0
\q\mbox{ in }\; D(A)'.
\end{equation}
Since $\f_j\in D(A)$, from
\eqref{11.2-eq2}, we find that
\vspace{-0.2cm}
\begin{equation}\label{11.2-eq1}
-\hat c_{j} (a+\mu_{j}) \int_0^1 \f_j(x)
\phi_{j}(x)dx-\sum_{k\in\dbZ^+\setminus\{j\}}
\hat c_k (a+\mu_{k}) \int_0^1 \f_j(x)
\phi_k(x)dx=\f_{j}'(0).
\end{equation}
For $j\neq k$, by \eqref{1.17-eq5},  \eqref{defhatb} and
\eqref{11.7-eq1}, we have that
$$
\begin{array}{ll}\ds
\q\hat c_k (a+\mu_{k})
\int_0^1 \f_j(x) \phi_k(x)dx\\
\ns\ds =\hat c_k (a+\mu_{k})
\int_0^1 \f_j(x)\big\{ \f_{k}(x) +(-A+\mu_{k}+a)^{-1}\big[a_{k}\hat b(x)-a\f_{k}(x) \,\big]\big\} dx\\
\ns\ds = \hat c_k (a+\mu_{k}) \int_0^1
\f_j(x)(-A+\mu_{k}+a)^{-1}a_{k} \f_j'(0)
 \f_j(x) dx\\
\ns\ds =  \hat c_k
(a+\mu_{k})(-\mu_j+\mu_{k}+a)^{-1}a_{k} \f_j'(0) \\
\ns\ds = \hat c_k
a\frac{a+\mu_{k}}{-\mu_j+\mu_{k}+a}\frac{\f'_{j}(0)}{\f'_{k}(0)}.
\end{array}
$$
This, together with \eqref{11.2-eq1},
implies that \vspace{-0.2cm}
\begin{equation}\label{11.3-eq3}
-\hat c_{j} (a+\mu_{j})  -
\sum_{k\in\dbZ^+\setminus\{j\}}\hat c_k a
\frac{a+\mu_{k}}{-\mu_j+\mu_k+a}\frac{\f'_{j}(0)}{\f'_{k}(0)}=\f'_{j}(0)\;
\mbox{ for all } j\in\dbZ^+.
\end{equation}
For $j,k\in\dbZ^+$, let \vspace{-0.2cm}
\begin{equation}\label{11.8-eq4}
 c_j \= -\frac{\hat
c_{j}(a+\mu_{j})}{\f'_{j}(0)}, \q a_{jk} \=
\frac{1}{-\mu_j+ \mu_{k}+a}.
\end{equation}
From \eqref{11.3-eq3} and \eqref{11.8-eq4},
we have that
\begin{equation}\label{11.4-eq1}
c_{j}+a\sum_{k\in\dbZ^+\setminus\{j\}} a_{jk}
c_k =1.
\end{equation}
Let us now estimate $c_j$.   From
\eqref{1.22-eq5}, \eqref{11.8-eq4} and
\eqref{11.4-eq1}, for every $j\in \dbZ^+$, we
have
\begin{equation}\label{11.7-eq15}
\begin{array}{ll}\ds
\Big|\sum_{k\in\dbZ^+\setminus\{j\}} a_{jk}
c_k\Big| \3n&\ds=
\Big|\sum_{k\in\dbZ^+\setminus\{j\}}\hat c_k
\frac{a+\mu_{k}}{-\mu_j+\mu_k+a}\frac{1}{\f'_{k}(0)}\Big|\\
\ns&\ds \leq
C\sum_{k\in\dbZ^+\setminus\{j\}}\Big|\hat c_k
\frac{k^3}{j^4 -
k^4}\Big|\\
\ns&\ds \leq
C\(\sum_{k\in\dbZ^+\setminus\{j\}}|\hat
c_k|^2\)^\frac{1}{2}\(\sum_{k\in\dbZ^+\setminus\{j\}}\Big|
\frac{k^3}{j^4 - k^4}\Big|^2\)^\frac{1}{2}.
\end{array}
\end{equation}
For any $j\in\dbZ^+$,
\begin{equation}\label{11.7-eq16}
\sum_{k\in\dbZ^+\setminus\{j\}}\Big|
\frac{k^3}{j^4 - k^4}\Big|^2 \leq \sum_{k>
2j}\Big| \frac{k^3}{j^4 - k^4}\Big|^2 + \sum_{j<
k\leq 2j}\Big| \frac{k^3}{j^4 - k^4}\Big|^2  +
\sum_{0< k < j}\Big| \frac{k^3}{j^4 -
k^4}\Big|^2.
\end{equation}
We now estimate the three terms in the right
hand side of \eqref{11.7-eq16}. Firstly, we
have that
\begin{equation}\label{11.7-eq17}
\sum_{k>2j}\Big| \frac{k^3}{j^4 - k^4}\Big|^2
\leq \frac{16}{15}\sum_{k>2j} \frac{1}{k^2}\leq
\frac{C}{j}.
\end{equation}
Secondly,
\begin{equation}\label{11.7-eq18}
\begin{array}{ll}\ds
\q
\Big|\frac{k^3}{j^4-k^4}\Big|^2\leq C \sum_{1\leq l\leq j} \Big|\frac{8j^3}{j^4-(j+l)^4}\Big|^2 \leq  C \sum_{1\leq l\leq j}
\Big|\frac{8j^3}{4j^3l + 4 j^2l^2 +
4jl^3+l^4}\Big|^2 \\
\ns\ds  \leq C \sum_{1\leq l\leq j}
\Big|\frac{j}{jl+l^2}\Big|^2
\leq C\sum_{1\leq
l\leq j}
\Big|\(\frac{1}{l}-\frac{1}{j+l}\)\Big|^2 \leq C\sum_{1\leq l\leq j} \(
\frac{1}{l^2}-\frac{1}{(j+l)^2}\)\leq C.
\end{array}
\end{equation}
Thirdly, with similar arguments as for \eqref{11.7-eq18}, we can obtain that
\begin{equation}\label{11.7-eq19}
\sum_{0< k< j}\Big| \frac{k^3}{j^4 - k^4}\Big|^2
\leq C.
\end{equation}

From \eqref{11.7-eq15} to \eqref{11.7-eq19},
we find that there is a constant $C>0$ such
that, for all positive integer $j$,
\begin{equation}\label{11.7-eq20.1}
\Big|\sum_{k\in\dbZ^+\setminus\{j\}}
a_{jk}c_k\Big|\leq C.
\end{equation}
Combining \eqref{11.4-eq1} and
\eqref{11.7-eq20.1}, we know that there is
a constant $C>0$ such that, for all
$j\in\dbZ^+$,
\begin{equation}\label{11.7-eq21}
|c_j|\leq C.
\end{equation}

We now estimate $|c_j|$ for large $j$. From
\eqref{11.8-eq4},
 we get  that
\begin{equation}\label{11.7-eq7}
\begin{array}{ll}\ds
\sum_{k\in\dbZ^+\setminus\{j\}}\!|a_{jk}|\!
 =\!\!\sum_{k\in\dbZ^+\setminus\{j\}}\!\frac{1}{|-\mu_j\!+\mu_k\!+\!a|}   \leq
\sum_{k\in\dbZ^+\setminus\{j\}}\frac{C}{|j^4-k^4|}\\
\ns\ds \leq
\sum_{k>2j}\frac{C}{\big|j^4-k^4\big|} +
\sum_{j< k\leq
2j}\frac{C}{\big|j^4-k^4\big|} +
\sum_{0<k<j}\frac{C}{\big|j^4-k^4\big|}.
\end{array}
\end{equation}
Let us estimate the terms in the last line
of \eqref{11.7-eq7} one by one. First,
\begin{equation}\label{11.7-eq8}
\sum_{k>2j}\frac{1}{|j^4-k^4|}\leq
\frac{16}{15} \sum_{k>2j}\frac{1}{k^4}\leq
\frac{16}{15}
\sum_{k>2j}\frac{1}{k^2}\frac{1}{k(k-1)}\leq
\frac{16}{15}\frac{1}{(2j+1)4j^2} \leq
\frac{2}{15j^3}.
\end{equation}
Second,\vspace{-0.4cm}
\begin{equation}\label{11.7-eq9}
\begin{array}{ll}\ds
\sum_{j< k\leq
2j}\frac{1}{|j^4-k^4|}\3n&\ds\leq \sum_{1<
l\leq j}\frac{1}{|j^4-(j+l)^4|}\leq \sum_{1<
l\leq
j}\frac{1}{4j^3l+6j^2l^2+4jl^3+l^4}\\
\ns&\ds \leq \frac{1}{4}\sum_{1< l\leq
j}\frac{1}{j^3l+j^2l^2}\leq
\frac{1}{4j^3}\sum_{1< l\leq
j}\(\frac{1}{j+l} + \frac{1}{l}\) \leq
\frac{\ln j}{2j^3}.
\end{array}
\end{equation}
Similarly, we can obtain that
\vspace{-0.3cm}
\begin{equation}\label{11.7-eq10}
\sum_{0<k<j}\frac{1}{\big|j^4-k^4\big|}\leq
\frac{\ln j }{2j^3}.
\end{equation}
From \eqref{11.7-eq7} to \eqref{11.7-eq10},
we know that there is a constant $C>0$ such that,
for all positive integer
$j$,\vspace{-0.2cm}
\begin{equation}\label{11.7-eq12}
\sum_{k\in\dbZ^+\setminus\{j\}}|a_{jk}|
\leq \frac{C(1+\ln j)}{j^3}.
\end{equation}
Combining \eqref{11.4-eq1}, \eqref{11.7-eq21}
and \eqref{11.7-eq12}, we obtain that
\vspace{-0.2cm}
\begin{equation}\label{11.7-eq22}
\left|c_j-1\right| \leq
\frac{C(1+\ln j)}{j^3}.
\end{equation}

Let us now turn to the construction of
$k(\cd,\cd)$. From \eqref{1.17-eq5},
\eqref{10.30-eq7} and \eqref{11.7-eq22}, one has
that \vspace{-0.2cm}
\begin{equation}\label{11.7-eq23}
\begin{array}{lcl}\ds
\q\vspace{-0.3cm}\sum_{j\in\dbZ^+}
\int_0^1\big|\f_{j}(x)-c_j\phi_{j}(x)\big|^2
dx\\
\ns\ds \vspace{-0.3cm}\leq \ds
2\sum_{j\in\dbZ^+}
\int_0^1\big|(1-c_j)\f_{j}(x)\big|^2 dx +
2\sum_{j\in\dbZ^+}
\int_0^1|c_j|^2\big|\f_{j}(x)-\phi_{j}(x)\big|^2
dx
\\
\ns\ds \vspace{-0.3cm} \leq
C\sum_{j\in\dbZ^+}\(\frac{\ln^2 j}{j^6}
+\frac{1}{j^6}\) <+\infty.
\end{array}
\end{equation}
Inequality \eqref {11.7-eq23} allows us to
define a $k(\cd,\cd)\in
L^2((0,1)\times(0,1))$ by \vspace{-0.1cm}
\begin{equation}
\label{defk} \vspace{-0.2cm}
k(x,y)\=\sum_{j\in\dbZ^+}[\f_{j}(x)-c_j\phi_{j}(x)]\f_j(y),\q
0\leq x,y\leq 1.
\end{equation}
Let us prove that $k(\cd,\cd)\in
H^2((0,1)\times (0,1))$.
Thanks to \eqref{1.17-eq5}, \eqref{10.30-eq2}
and \eqref{11.7-eq22}, we find that
\begin{equation}\label{11.7-eq24}
\begin{array}{ll}\ds
\q\int_0^1\int_0^1\Big|\sum_{j\in\dbZ^+}
[\f_{j}'(x)-c_j\phi'_{j}(x)]\f_j(y)\Big|^2
dxdy\\
\ns\ds  =  \sum_{j\in\dbZ^+}
\int_0^1\big|\f_{j}'(x)-c_j\phi_{j}'(x)\big|^2
dx \\
\ns\ds \leq\sum_{j\in\dbZ^+}
\int_0^1\big|(1-c_j)\f_{j}'(x)\big|^2 dx +
\sum_{j\in\dbZ^+}
\int_0^1|c_j|^2\big|\f_{j}'(x)-\phi_{j}'(x)\big|^2
dx
\\
\ns\ds \leq C\sum_{j\in\dbZ^+}\(\frac{\ln^2
j}{j^4} +\frac{1}{j^4}\) <+\infty,
\end{array}
\end{equation}
which implies that
\begin{equation}\label{kxL2}
k_x(\cdot,\cdot)\in L^2((0,1)\times(0,1)).
\end{equation}

Utilizing  \eqref{1.17-eq5}, \eqref{11.7-eq1}
and \eqref{11.7-eq21}, we get that
\vspace{-0.2cm}
\begin{equation}\label{11.7-eq25}
\begin{array}{l}\ds
\q\int_0^1\int_0^1\Big|\sum_{j\in\dbZ^+}
c_j[\f_{j}(x)-\phi_{j}(x)]\f_j'(y)\Big|^2dxdy
\\
[-1mm]\ns \ds= \int_0^1\sum_{j\in\dbZ^+}
j^2|c_j|^2|\f_{j}(x)-\phi_j(x)\big|^2dx \\
\\
[-4mm]\ns \ds =
\int_0^1\sum_{j\in\dbZ^+}\!\!j^2|c_j|^2
|(-A+\mu_{j}\!+\!a)^{-1}(a_{j}\hat b \!-\! a
\f_{j} )(x)|^2 dx  \\
\ns \ds = \int_0^1 \sum_{j\in\dbZ^+}j^2|c_j|^2
\Big|\sum_{l\in\dbZ^+\setminus\{j\}}(-\mu_l+\mu_{j}+a)^{-1}
\frac{a\f'_l(0)}{\f'_{j}(0)} \f_l(x)\Big|^2 dx
\\
\ns \ds =\sum_{l\in \dbZ^+}l^2
\sum_{l\in\dbZ^+\setminus\{l\}}|c_j|^2\left|(-\mu_l+\mu_{j}+a)^{-1}
\right|^2.
\end{array}
\end{equation}
Similarly to the proof of \eqref{11.7-eq12},
we can obtain that \vspace{-0.2cm}
\begin{equation}\label{11.7-eq26}
\sum_{j\in\dbZ^+\setminus\{l\}}\big|(-\mu_l+\mu_{j}+a)^{-1}
\big|^2 \leq \frac{C}{l^6},
\end{equation}
Combining \eqref{11.7-eq21}, \eqref{11.7-eq25}
and \eqref{11.7-eq26}, one has that
\vspace{-0.2cm}
\begin{equation}\label{11.7-eq27}
\begin{array}{ll}\ds
\int_0^1\int_0^1\Big|\sum_{j\in\dbZ^+}
c_j\big[\,\f_{j}(x)-\phi_{j}(x)\big]\f_j'(y)\Big|^2dxdy
\leq C \sum_{l\in\dbZ^+} \frac{1}{l^4}<+\infty.
\end{array}
\end{equation}
From  \eqref{11.7-eq22} and
\eqref{11.7-eq27}, it follows that
\vspace{-0.2cm}
\begin{equation}\label{11.8-eq1}
\begin{array}{ll}\ds
\q\int_0^1\int_0^1\Big|\sum_{j\in\dbZ^+}
\big[\,\f_{j}(x)-c_j\phi_{j}(x)\big]\f_j'(y)\Big|^2
dxdy\\
\ns\ds  \leq
2\int_0^1\!\int_0^1\!\Big|\!\sum_{j\in\dbZ^+}(1\!-\!c_j)\f_{j}(x)\f_j'(y)\Big|^2
dxdy \!+\!
2\int_0^1\!\int_0^1\!\Big|\!\sum_{j\in\dbZ^+}\!
c_j\big[\,\f_{j}(x)\!-\!\phi_{j}(x)\big]\f_j'(y)\Big|^2dxdy
\\
\ns\ds \leq 2\sum_{j\in\dbZ^+}
\int_0^1\big|(1-c_j)\f_{j}'(y)\big|^2 dy +
2\int_0^1\int_0^1\Big|\sum_{j\in\dbZ^+}
c_j\big[\,\f_{j}(x)-\phi_{j}(x)\big]\f_j'(y)\Big|^2dxdy
\\
\ns\ds \leq C\(\sum_{j\in\dbZ^+}\frac{1+
\ln^2j}{j^4}\)
 <+\infty.
\end{array}
\end{equation}
This, together with \eqref{11.7-eq23}  and
\eqref{11.7-eq24}, deduces  that \vspace{-0.2cm}
\begin{equation}
\label{kyL^2} k_y(\cd,\cd)\in
L^2((0,1)\times (0,1)),
\end{equation}
which, together with \eqref{kxL2}, shows
that
$k(\cd,\cd)\in H^1((0,1)\times
(0,1))$. Clearly, $k(\cd,\cd)=0$ on the
boundary of $(0,1)\times (0,1)$. Thus, we
conclude that
\begin{equation}\label{kinH10}
k(\cd,\cd)\in
H^1_0((0,1)\times (0,1)).
\end{equation}
Furthermore, by \eqref{1.17-eq5},
\eqref{10.30-eq2-double}  and \eqref{11.7-eq22},
we find that \vspace{-0.1cm}
\begin{equation}\label{11.7-eq24.1}
\begin{array}{ll}\ds
\q\int_0^1\int_0^1\Big|\sum_{j\in\dbZ^+}
[\f_{j}''(x)-c_j\phi''_{j}(x)]\f_j(y)\Big|^2
dxdy\\
\ns\ds \vspace{-0.1cm} =  \sum_{j\in\dbZ^+}
\int_0^1\big|\f_{j}''(x)-c_j\phi_{j}''(x)\big|^2
dx \\
\ns\ds \vspace{-0.1cm}\leq\sum_{j\in\dbZ^+}
\int_0^1\big|(1-c_j)\f_{j}''(x)\big|^2 dx +
\sum_{j\in\dbZ^+}
\int_0^1|c_j|^2\big|\f_{j}''(x)-\phi_{j}''(x)\big|^2
dx
\\
\ns\ds \leq C\sum_{j\in\dbZ^+}\(\frac{\ln^2
j}{j^2} +\frac{1}{j^2}\) <+\infty,
\end{array}
\end{equation}
which shows that
\begin{equation}\label{kxL2.1}
k_{xx}(\cdot,\cdot)\in
L^2((0,1)\times(0,1)).
\end{equation}
Utilizing \eqref{1.17-eq5},  \eqref{11.7-eq1}
and \eqref{11.7-eq21}, we get that
\begin{equation}\label{11.7-eq25.1}
\begin{array}{ll}\ds
\q\int_0^1\int_0^1\Big|\sum_{j\in\dbZ^+}
c_j[\f_{j}(x)-\phi_{j}(x)]\f_j''(y)\Big|^2dxdy
\\
\ns\ds  =\int_0^1 \sum_{j\dbZ^+}
j^4|c_j|^2\big|\f_{j}(x)-\phi_j(x)\big|^2dx \\
\ns\ds =
\int_0^1\sum_{j\in\dbZ^+}\!\!j^4|c_j|^2\left|
 (-A+\mu_{j}\!+\!a)^{-1}(a_{j}\hat b \!-\!
a \f_{j} )(x)\right|^2dx\\
\ns\ds = \int_0^1
\sum_{\in\dbZ^+}j^4|c_j|^2
\left|\sum_{l\in\dbZ^+\setminus\{j\}}(-\mu_l+\mu_{j}+a)^{-1}
\frac{a\f'_l(0)}{\f'_{j}(0)} \f_l(x)\right|^2dx\\
\ns\ds =a^2\!\!\sum_{j\in\dbZ^+}j^2|c_j|^2
\sum_{l\in\dbZ^+\setminus\{j\}}l^2\big|(-\mu_l+\mu_{j}+a)^{-1}\big|^2
\\
\ns\ds =a^2\!\!\sum_{k\in\dbZ^+}l^2
\sum_{j\in\dbZ^+\setminus\{k\}}j^2|c_j|^2\big|(-\mu_l+\mu_{j}+a)^{-1}\big|^2.
\end{array}
\end{equation}
Similar to the proof of \eqref{11.7-eq12}, we
can obtain that \vspace{-0.2cm}
\begin{equation}\label{11.7-eq26.1}
\sum_{j\in\dbZ^+\setminus\{l\}}j^2\big|(-\mu_l+\mu_{j}+a)^{-1}\big|^2
 \leq \frac{C}{l^4}.
\end{equation}
Combining \eqref{11.7-eq21},
\eqref{11.7-eq25.1} and \eqref{11.7-eq26.1}, one
has
\begin{equation}\label{11.7-eq27.1}
\begin{array}{ll}\ds
\int_0^1\int_0^1\Big|\sum_{j\in\dbZ^+}
c_j\big[\,\f_{j}(x)-\phi_{j}(x)\big]\f_j''(y)\Big|^2dxdy
\leq C \sum_{l\in\dbZ^+} \frac{1}{l^2}<+\infty.
\end{array}
\end{equation}
From \eqref{1.17-eq5}, \eqref{11.7-eq22} and
\eqref{11.7-eq27.1}, we get that
\begin{equation}\label{11.8-eq1.1}
\begin{array}{ll}\ds
\q\int_0^1\int_0^1\Big|\sum_{j\in\dbZ^+}
\big[\,\f_{j}(x)-c_j\phi_{j}(x)\big]\f_j''(y)\Big|^2
dxdy\\
\ns\ds  \leq
2\int_0^1\int_0^1\Big|\sum_{j\in\dbZ^+}(1-c_j)\f_{j}(x)\f_j''(y)\Big|^2
dxdy\\
\ns\ds \q +
2\int_0^1\int_0^1\Big|\sum_{j\in\dbZ^+}
c_j\big[\,\f_{j}(x)-\phi_{j}(x)\big]\f_j''(y)\Big|^2dxdy
\\
\ns\ds \leq 2\sum_{j\in\dbZ^+}
\int_0^1\big|(1-c_j)\f_{j}''(y)\big|^2 dy +
2\int_0^1\int_0^1\Big|\sum_{j\in\dbZ^+}
c_j\big[\,\f_{j}(x)-\phi_{j}(x)\big]\f_j''(y)\Big|^2dxdy
\\
\ns\ds \leq
C\sum_{j\in\dbZ^+}\frac{1+(\ln j)^2}{j^2}
 <+\infty.
\end{array}
\end{equation}
According to  \eqref{11.7-eq24.1} and
\eqref{11.8-eq1.1}, we see that \vspace{-0.1cm}
\begin{equation}
\label{kyyL^2} k_{yy}(\cd,\cd)\in
L^2((0,1)\times (0,1)),
\end{equation}
which, together with  \eqref{kinH10}, deduces
that $k(\cd,\cd)\in H^2((0,1)\times (0,1))$.

Define
\begin{equation}
\label{defkn} k^{(n)}(x,y) \= \sum_{0<j\leq
n}\big[\,\f_{j}(x)-c_j\phi_j(x)\big]\f_j(y)
\;\mbox{ in }(0,1)\times (0,1).
\end{equation}
From \eqref{1.17-eq5}, \eqref{phijH6} and
\eqref{defkn}, one has
\begin{equation}\label{10.13-eq4}
\left(x\in (0,1)\mapsto
k^{(n)}_{xx}(x,\cd)\in
L^2(0,1)\right)\text{ is in } C^0([0,1];
L^2(0,1)).
\end{equation}
For any $m,n\in\dbZ^+$, $m<n$, we have that
\begin{equation}\label{10.13-eq1}
\begin{array}{ll}\ds
\int_0^1\|\sum_{m< j\leq
n}\big[\,\f_{j}''(x)-c_j\phi_j''(x)\big]\f_{j}(y)\|^2dy\\
\ns\ds =\sum_{m<j\leq
n}\big|\f_{j}''(x)-c_j\phi_j''(x)\big|^2\\
\ns \ds \leq 2\sum_{m<j\leq
n}\big|(1-c_j)\f_{j}''(x)\big|^2 + 2\sum_{m<
j\leq
n}|c_j|^2\big|\f_{j}''(x)-\phi_j''(x)\big|^2.
\end{array}
\end{equation}
By means of  \eqref{11.7-eq22}, we find
that
\begin{equation}\label{11.8-eq2}
\max_{x\in [0,1]}\sum_{m< |j|\leq
n}\big|(1-c_j)\f_{j}''(x)\big|^2 \leq
C\sum_{m< |j|\leq n}\frac{1+ (\ln j)^2|}{j^2}.
\end{equation}
From \eqref{11.7-eq1}, similarly to the
proof of \eqref{11.7-eq12}, we obtain that
\begin{equation}\label{11.8-eq3}
\begin{array}{ll}\ds
\q\max_{x\in[0,1]}\sum_{m<j\leq
n}|c_j|^2\big|\f_{j}''(x)-\phi_j''(x)\big|^2\\
\ns\ds = \max_{x\in[0,1]}\sum_{m<j\leq
n}|c_j|^2\Big|\sum_{l\in\dbZ^+\setminus\{j\}}(-\mu_l+\mu_{j}+a)^{-1}
\frac{a\f'_l(0)}{\f'_{j}(0)}\f_l''(x)\Big|^2\\
\ns\ds \leq C \sum_{m<j\leq
n}\sum_{l\in\dbZ^+\setminus\{j\}}\Big|\frac{l^3}{(j^4-l^4)j}\Big|^2
\\
\ns\ds \leq C\sum_{m<j\leq n}\frac{\ln^2
j}{j^2}.
\end{array}
\end{equation}
Combining \eqref{11.8-eq2} and
\eqref{11.8-eq3}, we get that
\begin{equation}\label{10.13-eq5}
\left\{ x\in (0,1) \mapsto
k^{(n)}_{xx}(x,\cd)\in
L^2(0,1)\right\}_{n=1}^{+\infty} \mbox{ is
a Cauchy sequence in }C^0([0,1];L^2(0,1)),
\end{equation}
which deduces that \eqref{regkx} holds.
Proceeding as in the proofs of \eqref{kyyL^2}
and of \eqref{10.13-eq5}, one gets that
\begin{equation}\label{cvkny}
\left\{ y\in (0,1) \mapsto
k^{(n)}_{yy}(\cd,y)\in
L^2(0,1)\right\}_{n=1}^{+\infty} \mbox{ is a
Cauchy sequence in }C^0([0,1];L^2(0,1)),
\end{equation}
which also gives \eqref{regky}. Moreover,
\eqref{1.17-eq5}, \eqref{defk}, \eqref{defkn}
and \eqref{cvkny} imply that
\begin{equation}\label{11.8-eq6}
k_{yy}(\cd,0)=\lim_{n\to+\infty}k_{yy}^{(n)}(\cd,0)=0
\;\mbox{ in } L^2(0,1).
\end{equation}
Similarly, one can show that
\begin{equation}\label{11.8-eq7}
k_{yy}(\cd,1)=\lim_{n\to+\infty}k_{yy}^{(n)}(\cd,1)=0
\;\mbox{ in } L^2(0,1).
\end{equation}
Let us point that
\begin{equation}\label{ky0=0}
k_y(\cdot,0)=0 \text{ in } L^2(0,L).
\end{equation}
Indeed, from \eqref{defhatb}, \eqref{11.7-eq1},
\eqref{11.8-eq4}, \eqref{11.4-eq1},
\eqref{11.7-eq21}, \eqref{11.7-eq12} and
\eqref{defk}, one has, in $D(A)'$,
\begin{equation*}
\begin{array}{ll}
k_y(x,0)\3n&=\ds
\sum_{j\in\dbZ^+}[\f_{j}(x)-c_j\phi_{j}(x)]\f_j'(0)
\\
\ns&=\ds \sum_{j\in\dbZ^+}(1-c_j)\f_j'(0)\f_j(x)
-
\sum_{j\in\dbZ^+}\sum_{k\in\dbZ^+\setminus\{j\}}c_ja_{kj}a
\f_{k}'(0)\f_{k}(x)
\\
\ns&=\ds \sum_{j\in\dbZ^+}(1-c_j)\f_j'(0)\f_j(x)
- \sum_{k\in\dbZ^+}(1-c_k)\f_k'(0)\f_k(x)
\\
\ns &=\ds 0.
\end{array}
\end{equation*}

Let us now prove that $k(\cd,\cd)$ satisfies
\eqref{deftransposition}. From \eqref{1.22-eq5}
and \eqref{1.17-eq5} one has
\begin{equation}\label{11.11-eq1}
\begin{array}{ll}\ds
\q(\pa_{xxxx}+\l\pa_{xx} -
\pa_{yyyy}-\l\pa_{yy} +
a)(1-c_j) \f_j(x) \f_j(y)\\
\ns\ds =(-\mu_j+\mu_j+a) (1-c_j)\f_j(x)
\f_j(y)=a (1-c_j)\f_j(x) \f_j(y).
\end{array}
\end{equation}
 From \eqref{defhatb}, \eqref{11.8-eq5} and \eqref{11.7-eq1}, one
has
\begin{equation}\label{11.12-eq1}
\begin{array}{ll}\ds
\q(\pa_{xxxx}+\l\pa_{xx} -
\pa_{yyyy}-\l\pa_{yy} +
a)c_j[\f_j(x) -\phi_j(x)]\f_j(y)\\
\ns\ds = -(\pa_{xxxx}+\l\pa_{xx} -
\pa_{yyyy}-\l\pa_{yy} +
a)\sum_{k\in\dbZ^+\setminus\{j\}}c_j
\frac{a\f_k'(0)\f_k(x)}{\f_j'(0)(-\mu_k+\mu_j+a)}\f_j(y)\\
\ns\ds =
-c_ja\sum_{k\in\dbZ^+\setminus\{j\}}
\frac{\f_k'(0)}{\f_j'(0)} \f_k(x)\f_j(y)\\
\ns\ds =
c_ja\f_j(x)\f_j(y)-\frac{c_ja}{\f_j'(0)}
\sum_{k\in\dbZ^+}
\f_k'(0)\f_k(x)\f_j(y)\\
\ns\ds = c_ja \f_j(x)
\f_j(y)+\frac{c_ja}{\f_j'(0)}\d'_{x=0}\otimes\f_j(y)\q
\mbox{ in }D(A)'\otimes L^2(0,1).
\end{array}
\end{equation}
From \eqref{defkn}, \eqref{11.11-eq1} and
\eqref{11.12-eq1}, we get, in $D(A)'\otimes
L^2(0,1)$,
\begin{equation}
\label{eqkn}
\begin{array}{ll}\ds
k^{(n)}_{xxxx}(x,y)+\l
k^{(n)}_{xx}(x,y)-k^{(n)}_{yyyy}(x,y)-\l
k^{(n)}_{yy}(x,y)+a
k^{(n)}(x,y)\\
\ns\ds -\d'_{x=0}\otimes\sum_{0< j\leq
n}\frac{c_ja}{\f_j'(0)}\f_j(y)= a\sum_{0<j\leq
n} \f_j(x) \f_j(y).
\end{array}
\end{equation}
 From \eqref{def cE}, \eqref{1.17-eq5}, \eqref{defkn},
\eqref{eqkn}, one sees, using integrations by parts that, for any $\rho\in\cE\subset
D(A)\otimes L^2(0,1)$,
\begin{equation}\label{equationkn-testrho}
\begin{array}{ll}\ds
0\3n & \ds= \int_0^1 \int_0^1  \big[
\rho_{xxxx}(x,y) + \l\rho_{xx}(x,y) -
\rho_{yyyy}(x,y)- \l \rho_{yy}(x,y) + a\rho(x,y)
\big]k^{(n)}(x,y)dxdy\\
\ns&\ds \q  + \int_0^1
k_{y}^{(n)}(x,0)\rho_{yy}(x,0)dx - a
\sum_{0<j\leq n}\int_0^1 \int_0^1
\left(\sum_{0<j\leq n} \rho(x,y) \f_j(x)
\f_j(y)\right)dxdy.
\end{array}
\end{equation}
By \eqref{1.17-eq5}, \eqref{11.8-eq6}, \eqref{11.8-eq7} and
letting $n$ tends to $+\infty$ in
\eqref{equationkn-testrho}, we obtain that
\begin{equation}\label{7.30-eq1-new}
\begin{array}{ll}\ds
\int_0^1 \int_0^1  \big[ \rho_{xxxx}(x,y) +
\l\rho_{xx}(x,y) - \rho_{yyyy}(x,y) -
\l\rho_{yy}(x,y) + a\rho(x,y)
\big]k(x,y)dxdy\\
\ns \ds +  \int_0^1 k_{y}(x,0)\rho_{yy}(x,0)dx -
a \int_0^1 \rho(y,y) dy=0.
\end{array}
\end{equation}
which, together with \eqref{ky0=0}, gives
\begin{equation}\label{7.30-eq1-new1}
\begin{array}{ll}\ds
\int_0^1 \int_0^1  \big[ \rho_{xxxx}(x,y) +
\l\rho_{xx}(x,y) - \rho_{yyyy}(x,y) -
\l\rho_{yy}(x,y) + a\rho(x,y)
\big]k(x,y)dxdy\\
\ns \ds  - a\int_0^1 \rho(y,y) dy =0.
\end{array}
\end{equation}
Hence $k$ is a solution of \eqref{system3}. It
only remains to check  that the function $k$ is
real valued. This follows from the fact $\bar k$
is also a solution of
 \eqref{system3}. Hence by the uniqueness result of Step 1, we must have $\bar k=k$, which shows that
 the function $k$ is real valued. This concludes Step 2 and
therefore the proof of Theorem \ref{well-posed
lm}.
\endpf

%%%%%%%%%%%%%%%%%%%%%%%%%%%%%%%%%%%%%%%%

\section{Invertibility of $I-K$}
\label{secinvert}

%%%%%%%%%%%%%%%%%%%%%%%%%%%%%%%%%%%%%%%%

The goal of this section is to prove the
following theorem.
\begin{theorem}\label{invertible lm}
$I-K$ is an invertible operator from
$L^2(0,1)$ to $L^2(0,1)$.
\end{theorem}

{\it Proof of Theorem \ref{invertible lm}.}
Since $ k(\cd,\cd)\in L^2((0,1)\times (0,1))$,
we get that $K$ is a compact operator.
Furthermore, thanks to $k\in H^2((0,1)\times
(0,1))\cap H^{1}_0((0,1)\times (0,1))$, we know
that $K$ is a continuous linear map from
$L^2(0,1)$ into $H^2(0,1)\cap H^{1}_0(0,1)$.
Denote by $K^*$ the adjoint operator of $K$.
Then, it is easy to see that
$$
K^* (v)(x) = \int_0^1 k^*(x,y)v(y)dy
\q\mbox{ for any } v\in L^2(0,1),
$$
where $k^*$ is defined by
\begin{equation}\label{defk*}
k^*(x,y)\=k(y,x),\, (x,y)\in
(0,1)\times(0,1).
\end{equation}
From \eqref{system3} and \eqref{defk*}, we know
that $k^*(\cd,\cd)$ solves (in a sense which is
naturally adapted from \eqref{deftransposition})
\begin{equation}\label{system3.31}
\left\{
\begin{array}{ll}\ds
k^*_{yyyy} + \l k^*_{yy}  - k^*_{xxxx} -\l
k^*_{xx} +a k^* =a\d(y-x) &\mbox{ in }
(0,1)\times (0,1),\\
\ns\ds  k^*(0,y) =  k^*(1,y) =0 &\mbox{ on } (0,1),\\
\ns\ds  k^*_{xx}(0,y) =  k^*_{xx}(1,y) =0 &\mbox{ on } (0,1),\\
\ns\ds  k^*_{x}(0,y)  =0 &\mbox{ on } (0,1),\\
\ns\ds  k^*(x,0) =  k^*(x,1)=0 &\mbox{ on }
(0,1),\\
\ns\ds  k^*_{yy}(x,1) =0 &\mbox{ on }
(0,1).
\end{array}
\right.
\end{equation}
Furthermore, as a result of \eqref{regkx},
\eqref{regky} and \eqref{defk*}, we have the
following regularity for $k^*(\cd,\cd)$:
\begin{gather}
\label{9.5-eq5-1*} y\in (0,1)\mapsto
k_{yy}^*(\cd,y)\in L^2(0,1) \text{ is in }
C^0([0,1]; L^2(0,1)),
\\
\label{9.5-eq5-2*} x\in (0,1)\mapsto
k_{xx}^*(x,\cdot)\in L^2(0,1) \text{ is in }
C^0([0,1]; L^2(0,1)).
\end{gather}
From \eqref{defk*}, \eqref{system3.31} and
\eqref{9.5-eq5-2*}, we know that
\vspace{-0.15cm}
\begin{gather}\label{9.5-eq2}
 v\in C^2([0,1]) \text{ and }
v(0)=v_x(0)=v_{xx}(0)=v(1)=v_{xx}(1)=0, \q\forall v\in K^*(L^2(0,1)),
\\
\label{K*L2H2} K^* \text{ is a continuous linear
map from $L^2(0,1)$ to $H^2(0,1)$.}
\end{gather}

We claim that the spectral radius $r(K^*)$
of $K^*$ is equal to $0$. Otherwise, since $K^*:
L^2(0,1)\to L^2(0,1)$ is, as $K$, a compact
linear operator, it has a nonzero
eigenvalue $\a$. Then, there exists a
positive integer $n_0$ such that
\begin{equation}
\label{n0+1} \Ker(K^* - \a
I)^{n_0+1}=\Ker(K^* - \a I)^{n_0}.
\end{equation}
Let $ \cF\=\Ker(K^* - \a I)^{n_0}$. Since $K^*$
is a compact operator and $\alpha \not =0$,
$\cF$ is a finite dimensional vector space.
Moreover, since $\a\not =0$, one has $\cF\subset
K^*(L^2(0,1))$, which, together with
\eqref{9.5-eq2}, implies that
\begin{equation}
\label{vsurbord} \cF\subset\{v\in
C^{2}([0,1]):\,
v(0)=v_x(0)=v_{xx}(0)=v(1)=v_{xx}(1)=0\}.
\end{equation}
 Using the fact that, as $k$, $k^*\in
H^1_0((0,1)\times(0,1))$, we get that
\begin{equation}
\label{k*H-1} \text{$K^*$ can be extended
to be a continuous linear map from
$H^{-1}(0,1)$ into $L^2(0,1)$.}
\end{equation}
Denote by $\wt K^*$ this extension and
remark that, if $u\in H^{-1}(0,1)$ is such
that $\wt K^* u =\a u$, then $u\in
L^2(0,1)$. Thus,  $\Ker(\wt K^*-\a
I)\subset L^2(0,1)$ and $\Ker(\wt K^*-\a
I)=\Ker(K^*-\a I)$. Similarly, we have
\begin{equation}
\label{KertildeK*}
\left\{
\begin{array}{ll}
\ds\Ker(\wt K^*-\a I)^{n_0}\subset
L^2(0,1),\\
\ns\ds  \Ker(\wt K^*-\a I)^{n_0}=\cF, \\
\ns\ds\Ker(\wt K^*-\a I)^{n_0}=\Ker(\wt
K^*-\a I)^{n_0+1}.
\end{array}
\right.
\end{equation}
By \eqref{system3.31},
\begin{gather}
\label{commutation} K^*
(\pa_{xxxx}+\l\pa_{xx}) v =
(\pa_{xxxx}+\l\pa_{xx}) K^* v -a K^* v +
av, \, \forall\, v\in C^\infty_0(0,1).
\end{gather}
From \eqref{commutation} we get that $ K^*$ can
be extended to be a continuous linear map from
$H^{-4}(0,1)$ into $H^{-4}(0,1)$. Hence, by
interpolation and \eqref{K*L2H2}, we also get
that $ K^*$ can be extended as a continuous
linear map from $H^{-2}(0,1)$ into
$H^{-1}(0,1)$.   We denote by $\wh K^*$ this
extension. Then, from \eqref{KertildeK*}, one
has
\begin{equation}\label{9.5-eq4}
\left\{
\begin{array}{ll}
\ds\Ker(\wh K^*-\a I)^{n_0} =\Ker(K^*-\a
I)^{n_0} \subset
L^2(0,1),\\
\ns\ds  \Ker(\wh K^*-\a I)^{n_0}=\cF, \\
\ns\ds\Ker(\wh K^*-\a I)^{n_0}=\Ker(\wh
K^*-\a I)^{n_0+1}.
\end{array}
\right.
\end{equation}
Using a density argument,  \eqref{vsurbord} and
\eqref{commutation},
 we get, for every $v\in \cF$,
\begin{equation}\label{commutationwithA}
\wh K^* (\pa_{xxxx}+\l\pa_{xx}) v =
(\pa_{xxxx}+\l\pa_{xx}) \wh K^* v -a \wh K^* v +
a v  \text{ in }H^{-2}(0,1).
\end{equation}
 From  \eqref{vsurbord}, \eqref{commutationwithA},
 and induction on $n$, one
has, for every $v\in \cF$,
\begin{equation*}
(\wh K^*)^n (\pa_{xxxx}+\l\pa_{xx}) v=
(\pa_{xxxx}+\l\pa_{xx}) (\wh K^*)^n v - na (\wh
K^*)^n v + na(\wh K^*)^{n-1} v  \text{ in }
H^{-2}(0,1),
\end{equation*}
and therefore, for every polynomial $P$ and for every $v\in \cF$,
\begin{equation}\label{K*poly}
P(\wh K^*) (\pa_{xxxx}+\l\pa_{xx}) v =
(\pa_{xxxx}+\l\pa_{xx}) P(\wh K^*)v - a \wh K^*
P'(\wh K^*) v + a P'(\wh K^*) v  \text{ in }
H^{-2}(0,1).
\end{equation}
By virtue of \eqref{vsurbord}, \eqref{9.5-eq4} and
 \eqref{K*poly} with
$P(X)\=(X-\a)^{n_0+1}$, we see that
$(\pa_{xxxx}+\l\pa_{xx})\cF\subset \cF$.
Since $\cF$ is finite dimensional, we know
that $\pa_{xxxx}+\l\pa_{xx}$ has an
eigenfunction in $\cF$, that is, there
exist $\mu\in \dbC$ and
$\xi\in\cF\setminus\{0\}$ such that
$$
\left\{
\begin{array}{ll}\ds
(\pa_{xxxx}+\l\pa_{xx})\xi = \mu\xi &\mbox{ in } (0,1),\\
\ns\ds
\xi(0)=\xi_x(0)=\xi_{xx}(0)=\xi(1)=\xi_{xx}(1)=0.
\end{array}
\right.
$$
This leads to a contradiction with Corollary~\ref{cor-derivativenot0}. Then, we
know $r(K^*)=0$,  so $r(K)=0$. Hence,  the
real number $1$ belongs to the resolvent
set of $K$, which completes the proof of
Theorem~\ref{invertible lm}.
\endpf

\section{Proof of Theorem \ref{th1}}
\label{secproofth}

%%%%%%%%%%%%%%%%%%%%%%%%%%%%%%%%%%%%%%%%%%%%%%%%

This section is addressed to a proof of
Theorem \ref{th1}.

{\it Proof of Theorem \ref{th1}.} Let
$T>0$, which will be specified later on. Consider
the following equation

\vspace{-0.55cm}

\begin{equation}\label{12.22-eq6}
\left\{
\begin{array}{ll}\ds
v_{1,t} + v_{1,xxxx} + \l v_{1,xx}  +
v_1v_{1,x} = 0 &\mbox{ in }
[0,T]\times (0,1),\\
\ns\ds v_1(t,0) = v_1(t,1)=0 &\mbox{ on }
[0,T],\\
\ns\ds v_{1,xx}(t,0)=\int_0^1
k_{xx}(0,y)v_1(t,y)dy,\;v_{1,xx}(t,1)=0
&\mbox{ on } [0,T],
\\
\ns\ds v_1(0)=v^0 & \mbox{ in } (0,1).
\end{array}
\right.
\end{equation}
By Theorem \ref{well}, we know that there
exist  $r_T>0$ and $C_T>0$ such that, for all $v^0\in
L^2(0,1)$ with $|v^0|_{L^2(0,1)}\leq r_T$,
the equation \eqref{12.22-eq6} admits a
unique solution $v_1\in X_T$ and this solution satisfies

\vspace{-0.4cm}

\begin{equation}\label{12.22-eq2}
|v_1|_{X_T} \leq C_T|v^0|_{L^2(0,1)}.
\end{equation}

Let $w_1 = (I-K)v_1$. Then,  from the boundary
conditions of \eqref{system3} and
\eqref{12.22-eq6}, we have

\vspace{-0.7cm}

\begin{equation}\label{eq1}
\begin{array}{ll}\ds
w_{1,t}(t,x)\3n &\ds =v_{1,t}(t,x) -
\int_0^1
k(x,y)v_{1,t}(t,y)dy\\
\ns&\ds = v_{1,t}(t,x) + \int_0^1
k(x,y)\big[v_{1,yyyy}(t,y)+\l v_{1,yy}(t,y) + v_1(t,y)v_{1,y}(t,y)\big]dy \\
\ns&\ds = v_{1,t}(t,x) + \int_0^1
\big[k_{yyyy}(x,y)+\l k_{yy}(x,y)\big]
v_1(t,y)dy - \frac{1}{2}\int_0^1
k_y(x,y)v_1(t,y)^2dy,
\end{array}
\end{equation}
\vspace{-0.3cm}
\begin{equation}\label{eq2}
\l w_{1,xx}(t,x) =\l v_{1,xx}(t,x) -
\l\int_0^1 k_{xx}(x,y)v_1(t,y)dy
\end{equation}
and
\begin{equation}\label{eq3}
w_{1,xxxx}(t,x) = v_{1,xxxx}(t,x) - \int_0^1
k_{xxxx}(x,y)v_1(t,y)dy.
\end{equation}
 From \eqref{system3}, \eqref{12.22-eq6}, \eqref{eq1}, \eqref{eq2} and \eqref{eq3}, we obtain
\begin{equation}\label{eq4}
\begin{array}{ll}\ds
\q \!\!\!\!\!\vspace{-0.2cm} w_{1,t}(t,x)\!
+\! w_{1,xxxx}(t,x)\!\! +\!\! \l
w_{1,xx}(t,x)\!\! +\! a
w_1(t,x)\!+ \! \frac{1}{2}\!\int_0^1\!\! k_y(x,y)v_1(t,y)^2dy\!+\!v_1(t,x)v_{1,x}\!(t,x) \\
\ns\ds \vspace{-0.2cm} = a
v_1(t,x)-\int_0^1
v_1(t,y)\big[k_{xxxx}(x,y) +\l k_{xx}(x,y)
- k_{yyyy}(x,y)- \l k_{yy}(x,y) + a
k(x,y) \big] dy\\
\ns\ds = -\int_0^1
v_1(t,y)\big[k_{xxxx}(x,y)\! +\!\l
k_{xx}(x,y)\! -\! k_{yyyy}(x,y)\!-\! \l
k_{yy}(x,y)\! +\! a
k(x,y) \!-\!a \d(x-y) \big] dy \\
\ns\ds = 0.
\end{array}
\end{equation}
Hence, if we take the feedback control
$F(\cd)$ defined by \eqref{deffeedback},  we
get that $w_1$ solves
\begin{equation}\label{6.9-eq12}
\left\{ \begin{array}{ll}\ds w_{1,t} +
w_{1,xxxx} + \l
w_{1,xx} +  w_1w_{1,x} + a w_1 \\
\ns\ds   = - v_1v_{1,x} - \frac{1}{2}
\int_0^1 k_y(x,y)v_1(t,y)^2dy &\mbox{in }
[0,T]\times (0,1),\\
\ns\ds  w_1(t,0)=w_1(t,1)=0 &\mbox{on }[0,T],\\
\ns\ds  w_{1,xx}(t,0)=w_{1,xx}(t,1)=0 &\mbox{on
} [0,T].
\end{array}
\right.
\end{equation}
Multiplying the first equation of
\eqref{6.9-eq12} with $w_1$ and integrating on
$(0,1)$, we get, using integrations by parts,
\begin{equation}\label{6.9-eq13}
\begin{array}{ll}\ds
\q\frac{d}{dt}\int_0^1| w_1(t,x)|^2 dx \3n&=\ds
-2Q(w_1(t,\cdot))
 - \int_0^1
w_1(t,x)v_1(t,x)v_{1,x}(t,x)dx
\\
\ns&\displaystyle- \frac{1}{2}\int_0^1  w_1(t,x)
\[\int_0^1 k_y(x,y)v_1(t,y)^2dy\]dx,
\end{array}
\end{equation}
with
\begin{equation}\label{defQ}
Q(\varphi)\= \int_0^1|\varphi''(x)|^2 dx -
\l\int_0^1|\varphi'(x)|^2 dx +
a\int_0^1|\varphi(x)|^2 dx, \, \forall \varphi \in H^2(0,1)\cap H_0^1(0,1).
\end{equation}
Now we estimate the second and third terms in
the right hand side of \eqref{6.9-eq13}. First,
\begin{equation}\label{12.22-eq3}
\begin{array}{ll}\ds
\q\Big|\int_0^1
w_1(t,x)v_1(t,x)v_{1,x}(t,x)dx\Big|
\\
\ns\ds = \Big|\int_0^1
\[v_1(t,x)-\int_0^1 k(x,y)v_1(t,y)dy\]v_1(t,x)v_{1,x}(t,x)dx \Big|
\\
\ns\ds = \Big|\int_0^1 \(\int_0^1
k(x,y)v_1(t,y)dy\) v_1(t,x)v_{1,x}(t,x)dx
\Big|
\\
\ns\ds = \frac{1}{2}\Big|\int_0^1 \(\int_0^1
k_x(x,y)v_1(t,y)dy\) v_1(t,x)^2dx \Big|
\\
\ns\ds \leq
\frac{1}{2}\Big|\int_0^1|k_x(\cd,y)|^2dy\Big|_{L^\infty(0,1)}\(\int_0^1
v_1(t,x)^2dx\)^{\frac{3}{2}}
\\
\ns\ds \leq |(I-K)^{-1}|^3_{\cL(L^2(0,1))}
\frac{1}{2}\Big|\int_0^1|k_x(\cd,y)|^2dy\Big|_{L^\infty(0,1)}\(\int_0^1
w^2_1(t,x)dx\)^{\frac{3}{2}}.
\end{array}
\end{equation}
Next,
\begin{equation}\label{6.9-eq14}
\!\!\!\begin{array}{ll}\ds \q \Big|\int_0^1
w_1(t,x)\int_0^1
k_y(x,y)v_1(t,y)^2dy dx  \Big| \\
\ns\ds = \Big| \int_0^1 \[ v_1(t,x)-\int_0^1
k(x,y)v_1(t,y)dy \]\(\int_0^1
k_y(x,y)v_1(t,y)^2dy\)dx\Big|\\
\ns \ds \leq \!\Big|\!\int_0^1
\!\!\!v_1(t,x)\! \(\!\int_0^1\!\!
k_y(x,y)v_1(t,y)^2dy\!\)dx\Big| \!+\!
\Big|\! \int_0^1\!\! \(\!\int_0^1\!\!\!
k(x,y)v_1(t,y)dy \int_0^1\!\!
k_y(x,y)v_1(t,y)^2dy\! \)dx\Big|.
\end{array}
\end{equation}
The first term in the right hand side of
\eqref{6.9-eq14} satisfies that
\begin{equation}\label{10.15-eq1}
\begin{array}{ll}\ds
\q\Big| \int_0^1 v_1(t,x) \(\int_0^1
k_y(x,y)v_1(t,y)^2dy\)dx\Big|\\
\ns\ds = \Big| \int_0^1 v_1(t,y)^2
\(\int_0^1
k_y(x,y)v_1(t,x) dx\) dy \Big|\\
\ns\ds \leq \int_0^1 v_1(t,y)^2 \(\int_0^1
|k_y(x,y)|^2 dx\)^{\frac{1}{2}} \(\int_0^1
v_1(t,x)^2dx\)^{\frac{1}{2}}dy\\
\ns\ds \leq  \Big|\(\int_0^1 |k_y(x,\cd)|^2
dx\)^{\frac{1}{2}}\Big|_{L^\infty(0,1)}
\(\int_0^1 v_1(t,x)^2dx\)^{\frac{3}{2}}\\
\ns\ds \leq
|(I-K)^{-1}|_{\cL(L^2(0,1))}^{3}\Big|\(\int_0^1
|k_y(x,\cd)|^2
dx\)^{\frac{1}{2}}\Big|_{L^\infty(0,1)}
\(\int_0^1 w_1(t,x)^2dx\)^{\frac{3}{2}}.
\end{array}
\end{equation}
The second term in the right hand side of
\eqref{6.9-eq14} satisfies that
\begin{equation}\label{10.15-eq2}
\begin{array}{ll}\ds
\q\Big|\int_0^1 \(\int_0^1 k(x,z)v_1(t,z)dz
\int_0^1
k_y(x,y)v_1(t,y)^2dy\)dx\Big|\\
\ns\ds = \Big| \int_0^1 v_1(t,y)^2 \int_0^1
k_y(x,y)\(\int_0^1 k(x,z)v_1(t,z)dz\)dx dy \Big|\\
\ns\ds \leq \int_0^1 v_1(t,y)^2 \(\int_0^1
|k_y(x,y)|^2 dx\)^{\frac{1}{2}} \[\int_0^1
\(\int_0^1 k(x,z)v_1(t,z)dz\)^2dx\]^{\frac{1}{2}}dy\\
\ns\ds \leq  \Big|\(\int_0^1 |k_y(x,\cd)|^2
dx\)^{\frac{1}{2}}\Big|_{L^\infty(0,1)}\(\int_0^1\int_0^1
|k(x,y)|^2 dxdy\)^{\frac{1}{2}}
\(\int_0^1 v_1(t,x)^2dx\)^{\frac{3}{2}}\\
\ns\ds \leq
|(I\!-\!K)^{-1}|_{\cL(L^2(0,1))}^{3}\Big|\(\!\int_0^1\!
\!|k_y(x,\cd)|^2
dx\)^{\frac{1}{2}}\Big|_{L^\infty(0,1)}\!\(\!\int_0^1\!\!\int_0^1\!
|k(x,y)|^2 dxdy\)^{\frac{1}{2}}
\(\!\int_0^1\!\!\! w_1(t,x)^2dx\)^{\frac{3}{2}}.
\end{array}
\end{equation}

Let
$$
\begin{array}{ll}\ds
\wh C \3n&\ds=
\frac{1}{2}|(I-K)^{-1}|^3_{\cL(L^2(0,1))}
\Big|\int_0^1|k_x(\cd,y)|^2dy\Big|_{L^\infty(0,1)}
\\
\ns&\ds \q +
|(I-K)^{-1}|_{\cL(L^2(0,1))}^{3}\Big|\(\int_0^1
|k_y(x,\cd)|^2
dx\)^{\frac{1}{2}}\Big|_{L^\infty(0,1)}
\\
\ns&\ds \q +
|(I-K)^{-1}|_{\cL(L^2(0,1))}^{3}\Big|\(\int_0^1
|k_y(x,\cd)|^2
dx\)^{\frac{1}{2}}\Big|_{L^\infty(0,1)}\(\int_0^1\int_0^1
|k(x,y)|^2 dxdy\)^{\frac{1}{2}}.
\end{array}
$$
From  \eqref{6.9-eq13} to \eqref{10.15-eq2},
we get that
\begin{equation}\label{12.22-eq4}
\ds
\frac{d}{dt}\int_0^1 | w_1(t,x)|^2 dx
\leq -2Q(w_1(t,\cdot))
+\wh C\( \int_0^1
w_1(t,x)^2dx\)^{\frac{3}{2}}.
\end{equation}
By \eqref{nuassezpetit}, there exists $\nu'\in \mathbb{R}$ such that
\begin{equation}\label{propertynuprime}
\nu<\nu'\leq a+ j^4\pi^4 - \l j^2\pi^2, \q \forall j\in\dbZ^+.
\end{equation}
From \eqref{defQ} and the second inequality of \eqref{propertynuprime}, we have
\begin{equation}\label{Qassezgrand}
  Q(\varphi)\geq \nu' \int_0^1 |\varphi(x)|^2dx, \q \forall \varphi\in \in H^2(0,1)\cap H_0^1(0,1).
\end{equation}
 From \eqref{12.22-eq4} and \eqref{Qassezgrand}, one has
\begin{equation}\label{estevolw1}
\ds
\frac{d}{dt}\int_0^1 | w_1(t,x)|^2 dx
\leq -2\nu' \int_0^1 | w_1(t,x)|^2 dx
+\wh C\( \int_0^1
w_1(t,x)^2dx\)^{\frac{3}{2}}.
\end{equation}
By Theorem \ref{well} and the first inequality
of \eqref{propertynuprime}, there exists
$\d_1>0$ be such that, if $|w_1(0)|_{L^2(0,1)}
\leq \d_1$, then
\begin{equation*}
\wh C\( \int_0^1
w_1(t,x)^2dx\)^{\frac{1}{2}} \leq
2(\nu'-\nu),\quad \forall\, t\in [0,T].
\end{equation*}
This, together with \eqref{12.22-eq4} and
\eqref{estevolw1}, implies that
\begin{equation*}
\frac{d}{dt}\int_0^1 |w_1(t,x)|^2 dx \leq
-2\nu \int_0^1|w_1(t,x)|^2 dx,\quad \forall\,
t\in [0,T].
\end{equation*}
which gives us that
\begin{equation}\label{6.9-eq15}
|w_1(t,\cd)|_{L^2(0,1)} \leq e^{-\nu t}|w_1(0,\cd)|_{L^2(0,1)},
\q\forall\,t\in [0,T].
\end{equation}
Then, from Lemma \ref{invertible lm} and
\eqref{6.9-eq15}, we know  that,  if
$$
|v_1(0,\cd)|_{L^2(0,1)}\leq
\min\{|I-K|^{-1}_{\cL(L^2(0,1))}\d_1, r_T\},
$$
then
\begin{equation}\label{6.9-eq16}
|v_1(t,\cd)|_{L^2(0,1)} \leq
|(I-K)^{-1}|_{\cL(L^2(0,1))}|I-K|_{\cL(L^2(0,1))}e^{-\nu
t}|v_1(0,\cd)|_{L^2(0,1)},\q \forall\, t\in
[0,T].
\end{equation}
Now we choose $T>0$ such that $ e^{-\nu T} \leq
|(I-K)^{-1}|_{\cL(L^2(0,1))}|I-K|_{\cL(L^2(0,1))}$.
From \eqref{6.9-eq16}, we find that
$|v_1(T)|_{L^2(0,1)}\leq |v_1(0)|_{L^2(0,1)}
\leq r_T$. Thus, by Theorem \ref{well}, we know
that
\begin{equation}\label{12.22-eq5}
\left\{
\begin{array}{ll}\ds
v_{2,t} + v_{2,xxxx} + \l v_{2,xx}  +
v_2v_{2,x} = 0 &\mbox{ in }
[0,T]\times (0,1),\\
\ns\ds v_2(t,0) = v_2(t,1)=0 &\mbox{ on }
[0,T],\\
\ns\ds v_{2,xx}(t,0)=\int_0^1
k_{xx}(0,y)v_2(t,y)dy,\; v_{2,xx}(t,1)=0&\mbox{ on } [0,T], \\
\ns\ds v_2(0)= v_1(T) &\mbox{ in } (0,1),
\end{array}
\right.
\end{equation}
is well-posed. Furthermore, as for $v_1$ and
$w_1$, one can prove that $w_2\=(I-K)v_2$
satisfies
$$
|w_2(t,\cd)|_{L^2(0,1)} \leq e^{-\nu
t}|w_2(0,\cd)|_{L^2(0,1)}, \q\forall\, t\in
[0,T]
$$
and
$$
|v_2(T,\cd)|_{L^2(0,1)}\leq
|v_2(0,\cd)|_{L^2(0,1)}.
$$
Then, we can define $v_3$ and $w_3$ in a
similar manner. By induction, we can find
$v_n\in X_T$ ($n>1$), which solves
\begin{equation}\label{12.22-eq6.1}
\left\{
\begin{array}{ll}\ds
v_{n,t} + v_{n,xxxx} + \l v_{n,xx}  +
v_nv_{n,x} = 0 &\mbox{ in }
[0,T]\times (0,1),\\
\ns\ds v_n(t,0) = v_n(t,1)=0 &\mbox{ on }
[0,T],\\
\ns\ds v_{n,xx}(t,0)=\int_0^1
k_{xx}(0,y)v_n(t,y)dy,\; v_{n,xx}(t,1)=0 &\mbox{ on } [0,T], \\
\ns\ds v_n(0)= v_{n-1}(T) &\mbox{ in }
(0,1).
\end{array}
\right.
\end{equation}
Moreover, we have that $w_n = (I-K)v_n$
satisfies
\begin{equation}\label{12.22-eq7}
|w_n(t,\cd)|_{L^2(0,1)} \leq e^{-\nu
t}|w_{n}(0,\cd)|_{L^2(0,1)}= e^{-\nu
t}|w_{n-1}(T,\cd)|_{L^2(0,1)}
\end{equation}
and
$$
|v_n(T,\cd)|_{L^2(0,1)} \leq
|v_{n}(0,\cd)|_{L^2(0,1)}=|v_{n-1}(T,\cd)|_{L^2(0,1)}.
$$
Now we put
$$
v(t+(n-1)T,x) = v_n(t,x),\q
w(t+(n-1)T,x)=w_n(t,x) \q\mbox{ for }
(t,x)\in [0,T)\times [0,1].
$$
Then, it is an easy matter to see that $v$
solves \eqref{csystem1} and $w=(I-K)v$.
From \eqref{12.22-eq7}, we get that
\begin{equation}\label{12.22-eq8}
|w(t,\cd)|_{L^2(0,1)}\leq e^{-\nu
t}|w(0,\cd)|_{L^2(0,1)},\q\forall\, t\geq 0.
\end{equation}
This, together with $w=(I-K)v$, implies that
for all $t\geq 0$,
$$
|v(t,\cd)|_{L^2(0,1)}\leq e^{-\nu
t}|(I-K)^{-1}|_{\cL(L^2(0,1))}|I-K|_{\cL(L^2(0,1))}|v(0,\cd)|_{L^2(0,1)}\leq
Ce^{-\nu t}|v(0,\cd)|_{L^2(0,1)}.
$$
Let
$\d_0=\min\{|(I-K)|^{-1}_{\cL(L^2(0,1))}\d_1,
r_T\}$. Then, we know that for any $v^0\in
L^2(0,1)$ with $|v^0|_{L^2(0,1)}\leq \d_0$, the
equation \eqref{csystem1}  admits a solution $ v
\in X_T$. Furthermore, we have
\begin{equation}\label{12.22-eq9}
|v(t,\cd)|_{L^2(0,1)}\leq Ce^{-\nu
t}|v(0,\cd)|_{L^2(0,1)},\q\forall\, t\geq 0,
\end{equation}
which concludes the proof of Theorem~\ref{th1}. \endpf

\appendix
\section{Well-posedness of the Cauchy problem associated to the K-S control system}
\label{appendix-wellposed}

\q\;This section is devoted to a proof of
Theorem \ref{well}. Let $T>0$. We consider the
following linearized K-S equation with
non-homogeneous boundary condition:
\begin{equation}\label{12.19-eq1}
\left\{
\begin{array}{ll}\ds
u_t + u_{xxxx} + \l u_{xx} =\tilde h &\mbox{ in
}
[0,T]\times (0,1),\\
\ns\ds u(t,0)=u(t,1)=0 &\mbox{ on } [0,T],\\
\ns\ds u_{xx}(t,0) = h(t),\; u_{xx}(t,1)=0
&\mbox{
on } [0,T],\\
\ns\ds u(0)=u^0 &\mbox{ in }(0,1).
\end{array}
\right.
\end{equation}
Here $u^0\in L^2(0,1)$, $h\in L^2(0,T)$ and
$\tilde h\in L^2(0,T;L^2(0,1))$.

We first prove the following result:

\begin{lemma}\label{lm3}
Let $u^0\in L^2(0,1)$. There exists a unique
solution  $u\in C^0([0,T];L^2(0,1))\cap
L^2(0,T;H^2(0,1))$ of \eqref{12.19-eq1} such
that $u(0,\cdot)=u^0(\cdot)$. Moreover, there
exists $C_1\geq 1$, independent of  $\tilde h\in
L^1(0,T;L^2(0,1))$ and $u^0\in L^2(0,1)$, such
that
\begin{equation} \label{est1}
|u|_{X_T} \leq C_1\big(|u^0|_{L^2(0,1)} +
|h|_{L^2(0,T)} + |\tilde
h|_{L^2(0,T;L^2(0,1))}\big).
\end{equation}
\end{lemma}
{\it Proof of Lemma \ref{lm3}.} First, we assume
$u^0\in D(A)$ and $h\in H^1(0,T)$. Consider the
following equation:
\begin{equation}\label{1.23-eq2}
\left\{
\begin{array}{ll}\ds
\hat u_t + \hat u_{xxxx} + \l \hat u_{xx} = \hat
h &\mbox{ in }
[0,T]\times (0,1),\\
\ns\ds \hat u(t,0)=\hat u(t,1)=0 &\mbox{ on } [0,T],\\
\ns\ds \hat u_{xx}(t,0) = \hat u_{xx}(t,1)=0
&\mbox{ on }
[0,T],\\
\ns\ds \hat u(0)=\hat u^0,
\end{array}
\right.
\end{equation}
where
\begin{gather}\label{defhatu0}
\hat u^0\= u^0 + \frac{1}{6}(x^3 - 3x^2 +
2x)h(0),
\\
\label{defhath} \hat h = \tilde h +\lambda (x-1)
h(t) + \frac{1}{6}(x^3 - 3x^2 + 2x)h'(t).
\end{gather}
By the classical semigroup theory, we know that
\eqref{1.23-eq2} admits a unique solution $\hat
u$ in $C^0([0,T];D(A))$. Let
\begin{equation}\label{defuversushatu}
u(t,x) \= \hat u(t,x) -  \frac{1}{6}(x^3 - 3x^2
+ 2 x)h(t).
\end{equation}
Then it is easy to see that $u(\cd,\cd)$ solves
\eqref{12.19-eq1}.

Furthermore, multiplying both sides of
\eqref{12.19-eq1} by $u$ and integrating the
product in $(0,t)\times (0,1)$, we get that
\begin{equation}\label{1.23-eq4}
\begin{array}{ll}\ds
\int_0^1 u(t,x)^2dx - \int_0^1  u(0,x)^2dx +
2\int_0^t h(s)u_x(s,0)ds +
2\int_0^t\int_0^1 u_{xx}(s,x)^2 dxds\\
\ns\ds  +2\l\int_0^t\int_0^1 u_{xx}(s,x) u(s,x)
dxds - 2\l\int_0^t\int_0^1 u_x(s,x)^2 dxds =
2\int_0^t\int_0^1 \tilde h(s,x) u(s,x)dxds,
\end{array}
\end{equation}
which implies that, for every   positive real
number $\e$,
\begin{equation}\label{1.23-eq5}
\begin{array}{ll}\ds
\q\int_0^1 u(t,x)^2dx + 2\int_0^t\int_0^1
 u_{xx}(s,x)^2 dxds \\
\ns\ds \leq \int_0^1 u(0,x)^2dx +
2\int_0^t\int_0^1 \tilde h(s,x) u(s,x)dxds
-2\l\int_0^t\int_0^1 u_{xx}(s,x)
u(s,x) dxds \\
\ns\ds \q + \frac{1}{\e}\int_0^t h(s)^2ds +
\e\int_0^t u_x(s,0)^2ds.
\end{array}
\end{equation}
We have
\begin{equation}\label{2.11-eq1}
 -2\l\int_0^t\int_0^1 u_{xx}(s,x)
u(s,x) dxds \leq  \frac{1}{2}\int_0^t\int_0^1
 u_{xx}(s,x)^2 dxds +
2\l^2\int_0^t\int_0^1 u(s,x)^2 dxds.
\end{equation}
From \eqref{2.11-eq1} and integrations by parts,
we obtain that
\begin{equation}\label{2.11-eq2}
\begin{array}{lll}\ds
\q\int_0^t u_x(s,0)^2ds \3n&=\ds
\ds\int_0^t\int_0^1
\big[(x-1)u_x(s,x)^2\big]_x dxds \\
\ns&\ds = \ds2\int_0^t\int_0^1
(x-1)u_x(s,x)u_{xx}(s,x) dxds +
\int_0^t\int_0^1 u_x(s,x)^2 dxds\\
\ns\ds &\leq \ds  \int_0^t\int_0^1 u_{xx}(s,x)^2
dxds + 2\int_0^t\int_0^1 u_x(s,x)^2 dxds
\\
\ns\ds &\leq \ds  \int_0^t\int_0^1 u_{xx}(s,x)^2
dxds - 2\int_0^t\int_0^1 u(s,x)u_{xx}(s,x) dxds
\\
\ns\ds &\leq  \ds 2\int_0^t\int_0^1
u_{xx}(s,x)^2 dxds +  \int_0^t\int_0^1 u (s,x)^2
dxds.
\end{array}
\end{equation}
From \eqref{1.23-eq5} to \eqref{2.11-eq2}, by
choosing $\e= 1/4$, we find that
\begin{equation}\label{1.23-eq6}
\begin{array}{ll}\ds
\q\int_0^1  u(t,x)^2dx + \int_0^t\int_0^1
 u_{xx}(s,x)^2 dxds \\
\ns\ds \leq \int_0^1  u(0,x)^2dx + 4\int_0^t
 h(s)^2 ds + \int_0^t\int_0^1 \tilde
h(s,x)^2dxds + 2(1+\l^2)\int_0^t\int_0^1
u(s,x)^2dxds.
\end{array}
\end{equation}
This, together with Gronwall's inequality,
implies that
\begin{equation}\label{1.23-eq3}
|u|_{C^0([0,T];L^2(0,1))}\leq C\big( |
u^0|_{L^2(0,1)} + |h|_{L^2(0,T)} + |\tilde
h|_{L^2(0,T;L^2(0,1))} \big).
\end{equation}
According to \eqref{1.23-eq6} and
\eqref{1.23-eq3}, we get that
\begin{equation}\label{1.23-eq7.0}
|u|_{X_T}\leq C\big( |u^0|_{L^2(0,1)} +
|h|_{L^2(0,T)} + |\tilde h|_{L^2(0,T;L^2(0,1))}
\big).
\end{equation}
Now, by a standard density argument, we know
that for any $u^0\in L^2(0,1)$ and $h\in
L^2(0,T)$, \eqref{12.19-eq1} admits a solution
$u\in X_T$ which satisfies \eqref{1.23-eq7.0}.
The uniqueness of the solution follows from
\eqref{1.23-eq7.0}, which holds for every
solution of \eqref{12.19-eq1} in $X_T$. This
concludes the proof of Lemma~\ref{lm3}. \endpf

\begin{lemma}\label{lm4}
Let $z\in X_T$. Then $zz_x\in L^2(0,T;L^2(0,1))$
and the map $z\in X_T\mapsto zz_x\in
L^2(0,T;L^2(0,1))$ is continuous.
\end{lemma}
{\it Proof of Lemma~\ref{lm4}.} By the Sobolev
embedding theorem, we know that there is a
constant $\k>0$ such that
\begin{equation}\label{2.11-eq7}
|v|_{W^{1,\infty}(0,1))}\leq \k
|v|_{H^{2}(0,1)}.
\end{equation}
Then, we see that
\begin{equation}\label{1.23-eq9}
\begin{array}{ll}\ds
\q\int_0^T\int_0^1 |z(t,x)z_x(t,x)|^2 dxdt\3n
&\leq \ds \int_0^T
\(|z_x(t)|_{L^\infty(0,1)}^2\int_0^1
|z(t,x)|^2dx\)dt\\
\ns &\leq \ds \max_{t\in [0,T]}\int_0^1
|z(t,x)|^2dx\int_0^T
|z_x(t)|_{L^\infty(0,1)}^2 dt\\
\ns&=\ds |z|_{C^0([0,T];L^2(0,1))}^2
|z_x|_{L^2(0,T;L^\infty(0,1))}^2\\
\ns &\leq \ds C|z|_{C^0([0,T];L^2(0,1))}^2
|z|_{L^2(0,T;H^2(0,1))}^2.
\end{array}
\end{equation}
\endpf

{\it Proof of Theorem \ref{well}}\,: {\bf
Uniqueness of the solution.}

Assume that $u$ and $v$ are two solutions to
\eqref{csystem1}. Let $w=u-v$. Then we know that
$w$ solves
\begin{equation}\label{eqhatu-A}
\left\{
\begin{array}{ll}\ds
w_t  + w_{xxxx} + \l w_{xx}+ uw_x  +  wv_x =0
&\mbox{ in }
[0,T]\times (0,1),\\
\ns\ds w(t,0)=w(t,1)=0 &\mbox{ on } [0,T],\\
\ns\ds w_{xx}(t,0) = w_{xx}(t,1) =0 &\mbox{ on } [0,T],\\
\ns\ds w(0,x) = 0 &\mbox{ in } (0,1).
\end{array}
\right.
\end{equation}
Multiplying \eqref{eqhatu-A} by $w$, integrating
the product in $(0,t)\times (0,1)$ and
performing integration by parts, we get

\begin{equation}\label{2.11-eq3}
\begin{array}{ll}\ds
\q\int_0^1 w(t,x)^2dx  + 2\int_0^t\int_0^1
w_{xx}(s,x)^2 dxds
\\
\ns\ds =   -2\l  \int_0^t\int_0^1
w_{xx}(s,x)w(s,x) dxds + 2
\int_0^t\int_0^1 u(s,x)w_x(s,x)w(s,x) dxds \\
\ns\ds \q - 2\int_0^t\int_0^1
v(s,x)w_x(s,x)w(s,x)dxds.
\end{array}
\end{equation}
Since $w(s,\cd)\in H^2(0,1)\cap H_0^1(0,1)$, we
know that there is a constant $C_1>0$ such that
$|w(s)|_{H^2(0,1)}\leq
C_1|w_{xx}(s)|_{L^2(0,1)}$ for any $s\in [0,T]$.
Then, by Sobolev's embedding theorem, we get
that there is a constant $C_2>0$ such that
$|w(s)|_{W^{1,\infty}(0,1)}\leq
C_2|w_{xx}(s)|_{L^2(0,1)}$ for all $s\in [0,T]$.
Thus,
\begin{equation}\label{2.11-eq4}
\begin{array}{ll}\ds
\q\Big|\int_0^t\int_0^1 u(s,x)w_x(s,x)w(s,x)
dxds \Big|
\\
\ns\ds\leq
|u|_{C^0([0,T];L^2(0,1))}\int_0^t\(\int_0^1
|w_x(s,x)w(s,x)|^2 dx\)^{\frac{1}{2}}ds
\\
\ns\ds \leq
|u|_{C^0([0,T];L^2(0,1))}\int_0^t|w_x(s)|_{L^\infty(0,1)}\(\int_0^1
|w(s,x)|^2 dx\)^{\frac{1}{2}}ds
\\
\ns\ds \leq  \frac{\e}{4}
\int_0^t|w_x(s)|_{L^\infty(0,1)}^2 ds +
\frac{|u|_{C^0([0,T];L^2(0,1))}^2}{\e}
\int_0^t\int_0^1 |w(s,x)|^2 dxds
\\
\ns\ds \leq C_2\e \int_0^t\int_0^1
|w_{xx}(s,x)|^2 dxds +
\frac{|u|_{C^0([0,T];L^2(0,1))}^2}{\e}
\int_0^t\int_0^1 |w(s,x)|^2 dxds.
\end{array}
\end{equation}
Similarly, we have that
\begin{equation}\label{2.11-eq5}
\begin{array}{ll}\ds
\q\Big|\int_0^t\int_0^1 v(s,x)w_x(s,x)w(s,x)
dxds \Big|
\\
\ns\ds \leq C_2\e \int_0^t\int_0^1
|w_{xx}(s,x)|^2 dxds +
\frac{|v|_{C^0([0,T];L^2(0,1))}^2}{\e}
\int_0^t\int_0^1 |w(s,x)|^2 dxds.
\end{array}
\end{equation}
Moreover
\begin{equation}\label{2.11-eq6}
-2\l \int_0^t\int_0^1 w_{xx}(s,x)w(s,x) dxds
\leq \int_0^t\int_0^1 w_{xx}(s,x)^2 dxds + \l^2
\int_0^t\int_0^1 w(s,x)^2 dxds.
\end{equation}
From \eqref{2.11-eq3} to \eqref{2.11-eq6} and
taking $\e =C_2/4$, we get that, for any $t\in
[0,T]$,
\begin{equation}\label{1.4-eq3}
\int_0^1 |w(t,x)|^2 dx \leq \(\l^2 +
\frac{|u|_{C^0([0,T];L^2(0,1))}^2+|v|_{C^0([0,T];L^2(0,1))}^2}{\e}
\) \int_0^t\int_0^1 |w(s,x)|^2 dxds.
\end{equation}
This, together with the Gronwall inequality,
implies that $w=0$ in $[0,T]\times (0,1)$.

\vspace{0.4cm}

{\bf Existence of the solution}

 Let us extend
$\tilde h$ and $h$ to be functions on
$(0,+\infty)\times (0,1)$ and $(0,+\infty)$ by
setting them to be zero on $(T,+\infty)\times
(0,1)$ and $(T,+\infty)$, respectively. Denote
by $|\!|F|\!|$ the norm of the continuous linear
map $F:L^2(0,1)\mapsto \dbR$. Set
\begin{gather}
\label{defT1}
T_1\=\min\Big\{\frac{1}{2C_1^{2}|\!|F|\!|^{2}},T\Big\}.
\end{gather}
Let $u\in X_T$. We know that
$$
h(\cd)\=  F(u(\cd)) \in L^2(0,T_1).
$$
Hence, for $v^0$ given in $L^2(0,1)$, we can
define a map
$$
\cJ: X_{T_1}\to X_{T_1}
$$
by
\noindent $\cJ(u)=v$, where $v\in X_{T_1}$
solves \eqref{12.19-eq1} with  $h(\cd) =
F(u(\cd))$ and $v(0,\cdot)=v^0(\cdot)$.

For $\tilde u, \hat u\in X_{T_1}$, from
\eqref{est1} and \eqref{defT1} one has that
$$
\begin{array}{ll}\ds
\q|\cJ(\tilde u)-\cJ(\hat u)|_{X_{T_1}}
\\
\ns\ds \leq  C_1 \big|F(\tilde u(t,\cd)) -
F(\hat u(t,\cd)) \big|_{L^2(0,T_1)}
\\
\ns\ds \leq  C_1 |\!|F|\!|\(\int_0^{T_1}\int_0^1
|\tilde u(t,y)- \hat u(t,y)|^2
dydt\)^{\frac{1}{2}}
\\
\ns\ds \leq \sqrt{T_1} C_1 |\!|F|\!| |\tilde u-
\hat u|_{C^0([0,T_1];L^2(0,1))} \\
\ns\ds \leq \frac{1}{\sqrt{2}}|\tilde u- \hat
u|_{X_{T_1}}.
\end{array}
$$
Hence, we get that $\cJ(\cd)$ is a contractive
map. By the Banach fixed point theorem, we know
that $\cJ(\cd)$ has a unique fixed point $v_1$,
which is the solution to the following equation
\begin{equation}\label{12.19-eq2}
\left\{
\begin{array}{ll}\ds
v_{1,t} + v_{1,xxxx} + \l v_{1,xx} = \tilde h
&\mbox{ in }
[0,T_1]\times (0,1),\\
\ns\ds v_1(t,0)=v_1(t,1)=0 &\mbox{ on } [0,T_1],\\
\ns\ds v_{1,xx}(t,0) = F(v_1(t,\cd)),\;v_{1,xx}(t,1) = 0 &\mbox{ on } [0,T_1],\\
\ns\ds v_1(0,x)=v^0(x) &\mbox{ in } (0,1).
\end{array}
\right.
\end{equation}
Using $C_1\geq 1$, \eqref{est1}, \eqref{defT1}
and \eqref{12.19-eq2}, we find that
$$
\begin{array}{ll}\ds
|v_1|_{C^0([0,T_1];L^2(0,1))} \3n &\leq \ds
C_1\[|v^0|_{L^2(0,1)} + \( \int_0^{T_1}
\big|F(v(t,\cd))\big|^2 dt\)^{1/2} + |\tilde
h|_{L^2(0,T_1;L^2(0,1))}\]
\\
\ns&\ds\leq   C_1\(|v^0|_{L^2(0,1)} + \sqrt{T_1}
|\!|F|\!| |v_1|_{C^0([0,T_1];L^2(0,1))} +
|\tilde h|_{L^2(0,T_1;L^2(0,1))}\)
\\
\ns&\ds\leq  C_1\big(|v^0|_{L^2(0,1)} + |\tilde
h|_{L^2(0,T_1;L^2(0,1))}\big) +
\frac{1}{\sqrt{2}}
|v_1|_{C^0([0,T_1];L^2(0,1))},
\end{array}
$$
which implies that
\begin{equation}\label{12.19-eq3}
|v_1|_{C^0([0,T_1];L^2(0,1))} \leq
\left(2+\sqrt{2}\right)C_1\big(|v^0|_{L^2(0,1)}
+ |\tilde h|_{L^2(0,T_1;L^2(0,1))}\big).
\end{equation}
Thanks to \eqref{est1}, \eqref{defT1},
\eqref{12.19-eq2} and \eqref{12.19-eq3}, we
obtain that
\begin{equation}\label{estim-v1}
\begin{array}{ll}\ds
|v_1|_{X_T} \3n &\ds\leq C_1\[ |v^0|_{L^2(0,1)}
+ \(\int_0^{T_1}
\big|F(v(t,\cd))\big|^2dt\)^{1/2} + |\tilde
h|_{L^2(0,T_1;L^2(0,1))}\]
\\
\ns &\ds \leq C_1\big( |v^0|_{L^2(0,1)} +
|\tilde h|_{L^2(0,T_1;L^2(0,1))}\big) +
C_1T_1^{1/2} |\!|F|\!|
|v_1|_{C^0([0,T_1];L^2(0,1))}
\\
\ns&\ds \leq C_1\big( |v^0|_{L^2(0,1)} + |\tilde
h|_{L^2(0,T_1;L^2(0,1))}\big) +
\frac{1}{\sqrt{2}}|v_1|_{C^0([0,T_1];L^2(0,1))}
\\
\ns&\ds \leq 4C_1\big( |v^0|_{L^2(0,1)} +
|\tilde h|_{L^2(0,T_1;L^2(0,1))}\big).
\end{array}
\end{equation}
By a similar argument, we can prove that the
following equation
\begin{equation}\label{12.19-eq2.1}
\left\{
\begin{array}{ll}\ds
v_{2,t} + v_{2,xxxx} + \l v_{2,xx} = \tilde h
&\mbox{ in }
[0,T_1]\times (0,1),\\
\ns\ds v_2(t,0)=v_2(t,1)=0 &\mbox{ on } [0,T_1],\\
\ns\ds v_{2,xx}(t,0) = F(v_2(t,\cd)),\;v_{2,xx}(t,1) = 0 &\mbox{ on } [0,T_1],\\
\ns\ds v_2(T_1,x)=v_1(T_1,x) &\mbox{ in } (0,1)
\end{array}
\right.
\end{equation}
admits a unique solution in
$C^0([T_1,2T_1];L^2(0,1))\cap
L^2(T_1,2T_1;H^2(0,1)\cap H^1_0(0,1))$.
Furthermore, this solution satisfies that
\begin{equation}\label{estim-v2}
\begin{array}{ll}\ds
\q\left(|v_2|_{C^0([T_1,2T_1];L^2(0,1))}^2 +
|v_2|^2_{L^2(T_1,2T_1;H^2(0,1)\cap
H^1_0(0,1))}\right)^{1/2}\\
\ns\ds \leq
 4C_1\big(
|v_1(T_1)|_{L^2(0,1)} + |\tilde
h|_{L^2(T_1,2T_1;L^2(0,1))}\big).
\end{array}
\end{equation}
By induction, we know that for an integer $n\geq
2$, the following equation
\begin{equation}\label{eqvn-A}
\left\{
\begin{array}{ll}\ds
v_{n,t} + v_{n,xxxx} + \l v_{n,xx} = \tilde h
&\mbox{ in }
[(n-1)T_1,nT_1]\times (0,1),\\
\ns\ds v_n(t,0)=v_n(t,1)=0 &\mbox{ on } [(n-1)T_1,nT_1],\\
\ns\ds v_{n,xx}(t,0) = F(v_n(t,\cd)),\;v_{n,xx}(t,1) = 0  &\mbox{ on } [(n-1)T_1,nT_1],\\
\ns\ds v_n((n-1)T_1,x)=
v_{n-1}((n-1)T_1,x)&\mbox{ in } (0,1),
\end{array}
\right.
\end{equation}
admits a unique solution. Moreover, one has
\begin{equation}\label{estim-vn}
\begin{array}{ll}\ds
\q\left(|v_n|_{C^0([(n-1)T_1,nT_1];L^2(0,1))}^2
+ |v_n|^2_{L^2((n-1)T_1,nT_1;H^2(0,1)\cap
H^1_0(0,1))}\right)^{1/2}
\\
\ns\ds\leq 4C_1\big(
|v_{n-1}((n-1)T_1)|_{L^2(0,1)} + |\tilde
h|_{L^2((n-1)T_1,nT_1;L^2(0,1))}\big).
\end{array}
\end{equation}
Let $n_0\in \dbZ^+$ be such that $(n_0-1)T_1
\leq T < n_0T_1$. Let
$$
v(t,x)\=v_n(t,x) \mbox{ for } t\in
[(n-1)T_1,nT_1)\cap [0,T],\q n=1,2,\cds,n_0, \,
x\in (0,1).
$$
Then $v$ is the solution to
\begin{equation}\label{eqv-A}
\left\{
\begin{array}{ll}\ds
v_{t} + v_{xxxx} + \l v_{xx} = \tilde h &\mbox{
in }
[0,T]\times (0,1),\\
\ns\ds v(t,0)=v(t,1)=0 &\mbox{ on }
[0,T],\\
\ns\ds v_{xx}(t,1) = F(v(t,y)),\;v_{xx}(t,1) = 0  &\mbox{ on } [0,T],\\
\ns\ds v(0,x)= v^0(x)&\mbox{ in } (0,1).
\end{array}
\right.
\end{equation}
Moreover,  there is a constant $C=C(T)>0$
(independent of $v^0\in L^2(0,1)$ and of $\tilde
h\in L^2(0,T;L^2(0,1))$) such that
\begin{equation}\label{estim-lin-kdv}
|v|_{X_T}^2  \leq C(T)\big( |v^0|_{L^2(0,1)}^2 +
|\tilde h|_{L^2(0,T;L^2(0,1))}^2\big).
\end{equation}

From now on, we assume that $v^0\in L^2(0,1)$
satisfies
\begin{equation}\label{conditionu0}
|v^0|_{L^2(0,1)}^2 \leq
\frac{5}{36C(T)^2\kappa^2},
\end{equation}
where $\k$ is the constant given in
\eqref{2.11-eq7}.

Let
$$
\cB\= \Big\{u \in X_T: \, |u|^2_{X_T} \leq
\frac{1}{6C(T)\kappa^2} \Big\}.
$$
Then $\cB$ is a nonempty closed subset of the
Banach space $X_T$. Let us define a map $\cK$
from $\cB$ to $X_T$ as follows:
$$
\cK(z) = v, \mbox{ where $v$ is the solution to
\eqref{eqv-A} with } \tilde h = zz_x.
$$
From Lemma \ref{lm4} and the above
well-posedness result for \eqref{eqv-A}, we know
that $\cK$ is well-defined. From
\eqref{estim-lin-kdv} and \eqref{conditionu0},
we have that
$$
\begin{array}{ll}\ds
|\cK(z)|^2_{X_T} \3n&\ds \leq
C(T)\big(|v^0|_{L^2(0,1)}^2 +
|zz_x|_{L^2(0,T;L^2(0,1))}^2\big)
\\
\ns&\ds \leq C(T)\big(|v^0|_{L^2(0,1)}^2 +
|z|^2_{C^0([0,T];L^2(0,1))}|z_x|^2_{L^2(0,T;L^\infty(0,1))}\big)
\\
\ns&\ds \leq C(T)\big( |v^0|_{L^2(0,1)}^2 +
\kappa^2 |z|^4_{X_T}\big)
\\
\ns&\ds \leq \frac{1}{6C(T)\kappa^2}
\end{array}
$$
and
$$
\begin{array}{ll}\ds
\q|\cK(z_1)-\cK(z_2)|^2_{X_T}
\\
\ns\ds \leq C(T)
|z_1z_{1,x}-z_2z_{2,x}|^2_{L^2(0,T;L^2(0,1))}=
C(T) |z_1z_{1,x}-z_1z_{2,x} +
z_1z_{2,x}-z_2z_{2,x} |^2_{L^2(0,T;L^2(0,1))}
\\
\ns\ds \leq 2C(T)\kappa^2 \big(
|z_1|^2_{C^0([0,T];L^2(0,1))}|z_1-z_2|^2_{L^2(0,T;H^2(0,1))}
+
|z_1-z_2|^2_{C^0([0,T];L^2(0,1)))}|z_2|^2_{L^2(0,T;H^2(0,1))}
\big)
\\
\ns\ds \leq
2C(T)\kappa^2\big(|z_1|^2_{C^0([0,T];L^2(0,1))}+|z_2|^2_{L^2(0,T;H^2(0,1))}\big)
|z_1-z_2|^2_{X_T} \\
\ns\ds \leq \frac{2}{3}|z_1-z_2|^2_{X_T}.
\end{array}
$$
Therefore, we know that $\cK$ is from $\cB$ to
$\cB$ and is contractive. Then, by the Banach
fixed point theorem, there is a (unique) fixed
point $v$, which is the solution to
\eqref{csystem1} with the initial datum
$v(0,\cdot)=v^0(\cdot)$.
\endpf

\section{On the approximate controllability of the linear K-S control system \eqref{10.4-eq1}}
\label{sec-app-contrllability}

In this section we prove the following theorem.

\begin{theorem}\label{th app}
System \eqref{10.4-eq1} is approximately
controllable if and only if $\l\notin\cN$.
\end{theorem}

{\it Proof of Theorem \ref{th app}.} Consider
the following equation
\begin{equation}\label{3.2-eq3}
\left\{
\begin{array}{ll}\ds
z_t + z_{xxxx} + \l z_{xx} = 0 &\mbox{ in } (0,T)\times(0,1),\\
\ns\ds z(t,0)=z(t,1)=0 &\mbox{ in } (0,T),\\
\ns\ds z_{xx}(t,1)=0,\; z_{xx}(t,0)=0
&\mbox{ in } (0,T),\\
\ns\ds z_{x}(t,0)=0  &\mbox{ in } (0,T),\\
\ns\ds z(0)=z_0 &\mbox{ in } (0,1),
\end{array}
\right.
\end{equation}
where $z_0\in L^2(0,1)$. Let us recall that the
approximate controllability of \eqref{10.4-eq1},
is equivalent to the property that the only
solution to \eqref{3.2-eq3} is zero (see, e.g.
\cite[Theorem 2.43]{2007-Coron-book}).

The ``if" part: Let $z_{0,j}=\int_0^1
z_0(x)\sqrt{2}\sin (j\pi x)dx$ for $j\in\dbZ^+$.
The solution to \eqref{3.2-eq3} reads
$$
z(t,x)=\sum_{j\in\dbZ^+}z_{0,j}e^{(-j^4\pi^4 +
\l j^2\pi^2)t}\sqrt{2}\sin (j\pi x).
$$
Then,
$$
z_x(t,0)=\sum_{j\in\dbZ^+}z_{0,j}e^{(-j^4\pi^4 +
\l j^2\pi^2)t}\sqrt{2}j\pi=0.
$$
Clearly, $z_x(t,0)$ is analytic in $(0,T]$.
Since $z_x(t,0)=0$ in $(0,T]$, we know that
\begin{equation}\label{3.2-eq5}
z_x(t,0)=\sum_{j\in\dbZ^+}z_{0,j}e^{(-j^4\pi^4 +
\l j^2\pi^2)t}\sqrt{2}j\pi=0 \mbox{ in
}(0,+\infty).
\end{equation}
 Let
$j_0\in\dbZ^+$ such that $-j_0^4\pi^4 + \l
j_0^2\pi^2 = \max_{j\in\dbZ^+}\{-j^4\pi^4 + \l
j^2\pi^2\}$. From \eqref{3.2-eq4}, we have that
$j_0$ is unique. Let us multiply both sides of
\eqref{3.2-eq5} by $e^{(j_0^4\pi^4 - \l
j_0^2\pi^2)t}$ and let $t$ tends to $+\infty$.
Then, we find that $z_{0,j_0}=0$. Now let
$j_1\in\dbZ^+$ be such that $-j_1^4\pi^4 + \l
j_1^2\pi^2 = \max_{j\in(\dbZ^+\setminus
\{j_0\})}\{-j^4\pi^4 + \l j^2\pi^2\}$. From
\eqref{3.2-eq4}, we have that $j_1$ is unique.
Let us multiply both sides of \eqref{3.2-eq5} by
$e^{(j_1^4\pi^4 - \l j_1^2\pi^2)t}$ and let
$t\to+\infty$. Then, we find that $z_{0,j_1}=0$.
By induction, we can conclude that $z_{0,j}=0$
for all $j\in\dbZ^+$, which implies that $z=0$
in $(0,T)\times (0,1)$.

The ``only if" part: If $\l\in\cN$, then there
are
 $j_0,k_0\in\dbZ^+$ with $j_0\neq k_0$ such that $\l =j_0^2\pi^2 + k_0^2\pi^2$. Thus,
we get $e^{(-j_0^4\pi^4 + \l
j_0^2\pi^2)t}=e^{[-k_0^4\pi^4 + \l
k_0^2\pi^2]t}=e^{j_0^2k_0\pi^4 t}$. Let us
choose
$$z(t,x)=e^{(-j_0^4\pi^4 + \l
j_0^2\pi^2)t}\sin (j_0\pi x) - e^{[-k_0^4\pi^4 +
\l k_0^2\pi^2]t}\frac{j_0}{k_0}\sin \big(k_0\pi
x\big).$$ Then, clearly, $z$ solves
\eqref{3.2-eq3} and $z$ is not zero. This
concludes the proof of Theorem \ref{th app}.
\endpf
As a corollary of Theorem~\ref{th app}, one has
the following corollary (which can be also
checked directly)
\begin{corollary}
\label{cor-derivativenot0} Assume that $\lambda
\not \in \mathcal{N}$ and that $\varphi \in
C^4([0,1];\mathbb{C})$ is such that, for some
$\mu \in \mathbb{C}$,
\begin{equation}\label{eqvarphi}
\left\{
\begin{array}{l}
  \varphi''''+\lambda \varphi''=\mu \varphi \text{ in } (0,1),
  \\
  \varphi(0)=\varphi(1)=\varphi'(0)=\varphi''(0)=\varphi''(1)=0,
\end{array}
\right.
\end{equation}
then, $\varphi=0$.
\end{corollary}
{\it Proof of Corollary
\ref{cor-derivativenot0}.} It suffices to remark
that, if \eqref{eqvarphi} holds, then
$z(t,x)=e^{-\mu t} \varphi(x)$ is a solution of
\eqref{3.2-eq3}. \endpf

%%%%%%%%%%%%%%%%%%%%%%%%%%%%%%%%%%%%%%%%%%%%%%%%%%%%%%%%%%%%%%%%%%%%%%%%%%%%%%%%%%%%

\section*{Acknowledgement} Qi L\"{u} would like to thank the Laboratoire Jacques-Louis
Lions (LJLL) at Universit\'{e} Pierre et
Marie Curie in Paris for its hospitality.
This work was carried out while he was
visiting the LJLL. We appreciate Dr.
Eduardo Cerpa for his useful comments which
help us improve this paper.

\end{document}